\documentclass{article}
\usepackage{graphicx}
\usepackage{amsmath, amsthm, amssymb, amscd}
\usepackage{cite}
\def\url#1{{\tt #1}}
\def\href#1#2{{#2}}

\def\clap#1{\hbox to 0pt{\hss#1\hss}}

\def\Vec#1{\vec{\mathbf{#1}}}
\def\VecS#1{\hat{\boldsymbol{#1}}}

\def\beq{\begin{equation}}
\def\eeq{\end{equation}}

\begin{document}

\title{ 
Cavatappi 2.0: More of the same but better
}
\author{
Robert T. Jantzen\\ 
Department of Mathematics and Statistics\\
Villanova University
}
\maketitle

\begin{abstract}
Innocent musing on geodesics on the surface of helical pasta shapes leads to a single continuous 4-parameter family of surfaces invariant under at least a 1-parameter symmetry group and which contains as various limits spheres, tori, helical tubes, and cylinders, all useful for illustrating various aspects of geometry in a visualizable setting that are important in special and general relativity. In this family the most aesthetically pleasing surfaces come from screw-rotating a plane cross-sectional curve perpendicular to itself, i.e., orthogonal to the tangent vector to a helix. If we impose instead this orthogonality in the Lorentzian geometry of 3-dimensional Minkowski spacetime with a timelike helical ``central" world line representing a circular orbit, we can model the Fermi Born rigid model of the classical electron in such an orbit around the nucleus, and visualize the Fermi coordinate grid and its intersection with the world tube of the equatorial circle of the spherical surface of the electron (suppressing one spatial dimension). This leads what we might playfully term ``relativistic pasta." This is a useful 2-dimensional stationary spacetime with closed spatial slices on which to illustrate the slicing and threading splittings of general relativity relative to a Killing congruence, like stationary axisymmetric spacetimes including the rotating black hole family.
\end{abstract}

\section{Introduction}

They say one thing leads to another... so it goes in mathematical musing. From the sphere to the torus \cite{bobdgn,bobtorus} to the cavatappi corkscrew pasta surface  \cite{pastadesign,blog,NYTimes,bobcavatappo10}, a little more whimsical thought produces a new improved cavatappo 2.0 that comes from enlarging this single continuous family even more to allow it to encompass the cylinder as well, another important example. This gives us a mathematical playground to learn about interesting aspects of geometry that pop up in special and  general relativity and spacetime geometry \cite{mtw}.

In an effort to find interesting examples for a rarely offered differential geometry course in my department, the family of tori presented themselves as the obvious first class of simple surface shapes that one encounters after the plane, cylinder, cone and sphere, each of which has a 3-parameter isometry group of maximal symmetry. The torus has the advantage of being the first nontrivial example of a surface created from a circle which only admits a 1-parameter isometry group, like the larger family of surfaces of revolution to which it belongs. However, less familiar than rotational symmetry is the screw symmetry \cite{wikiscrew} that actually has practical applications even today, several millenia since their first use by Archimedes and others \cite{archimedes}. My interest was sparked by the Mathematica analyzed pasta shapes of George Legendre in his coffee table book \textit{Pasta by Design} \cite{pastadesign}, where the corkscrew pasta shape called by many names (cavatappi, cellentani, gobetti, etc.) offered itself up as a great candidate for geodesic analysis using the same techniques valid for surfaces of revolution. 

While surfaces of revolution like the cylinder or sphere or torus have natural orthogonal coordinate grids formed by their meridians and parallels, the generalized meridians and parallels of  corkscrew surfaces are no longer orthogonal and require the use of a new orthogonal frame for the tangent spaces obtained simply by completing the square on the meridian or parallel differentials in the surface metric in order to tame the geometry into a manageable one. This is good practice for dealing with rotating black hole spacetime metrics in general relativity, where such a technique is standard practice \cite{mtw,mfg}, and these surfaces even provide examples of the closely related synchronization defect and Sagnac effect that accompany spacetimes with Killing congruences having nonzero vorticity \cite{clock} like the trajectories of corkscrew motions.
However, the corkscrew surfaces in Euclidean space become very interesting examples even in special relativity when one uses the flat Lorentzian metric instead of the Euclidean one on $R^3$, describing the deformation of a spherical body in a circular orbit undergoing Lorentz contraction along the direction of motion, while allowing a concrete realization of the Thomas precession of an electron in the adjustment of the surface grid to remove the twisting due to the torsion of the Frenet-Serret frame of the central helix (called the Fermi-Walker angular velocity of the frame). This  actually provides a low-dimensional exact example of Fermi coordinates \cite{fermicircles}. All of these applications provide a gateway for interesting geometry that is a useful foundation for theoretical physics on the one hand, or simply for popular exposure to ideas from that field which catch the public's imagination on the other. Either way, it gives a traditional subject a breath of new life that cannot be a bad thing.

\begin{figure}[t]
\typeout{*** EPS figure cavatappo}
\vglue-4cm
\begin{center}
\includegraphics[scale=0.5]{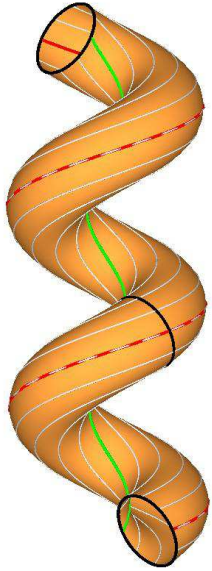}
\end{center}
\vglue-2cm

\caption{The new improved cavatappo/cellentano/gobbetto pasta surface (cavatappo 2.0), shown here in the orthogonally tilted case with parameter values: $a=3/2 b,c=4/5 b \to c/a=8/15$ and a central helix inclination angle of $\arctan(8/15)\approx 28.1^o$. A selection of equally spaced parallels are shown in white, together with the special parallels:  the outer/inner equators (red/green) and the northern/southern helices (cyan). Three meridians are shown in black.
}\label{fig:cavatappo}
\end{figure}

\begin{figure}[p] 
\typeout{*** EPS figure torus1}
\vglue0cm
\begin{center}
\includegraphics[scale=0.3]{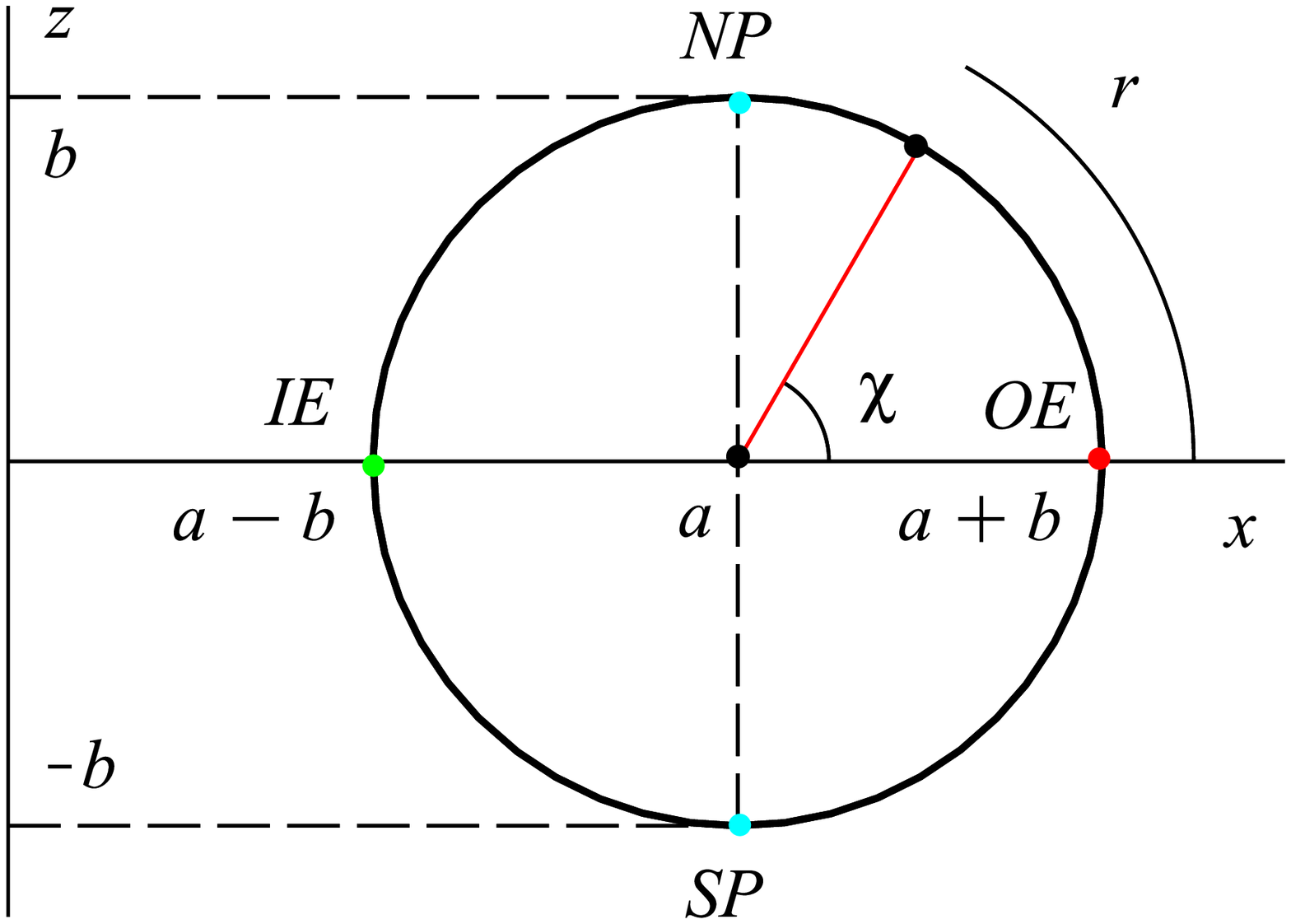}%\\[-1cm]
 \includegraphics[scale=0.3]{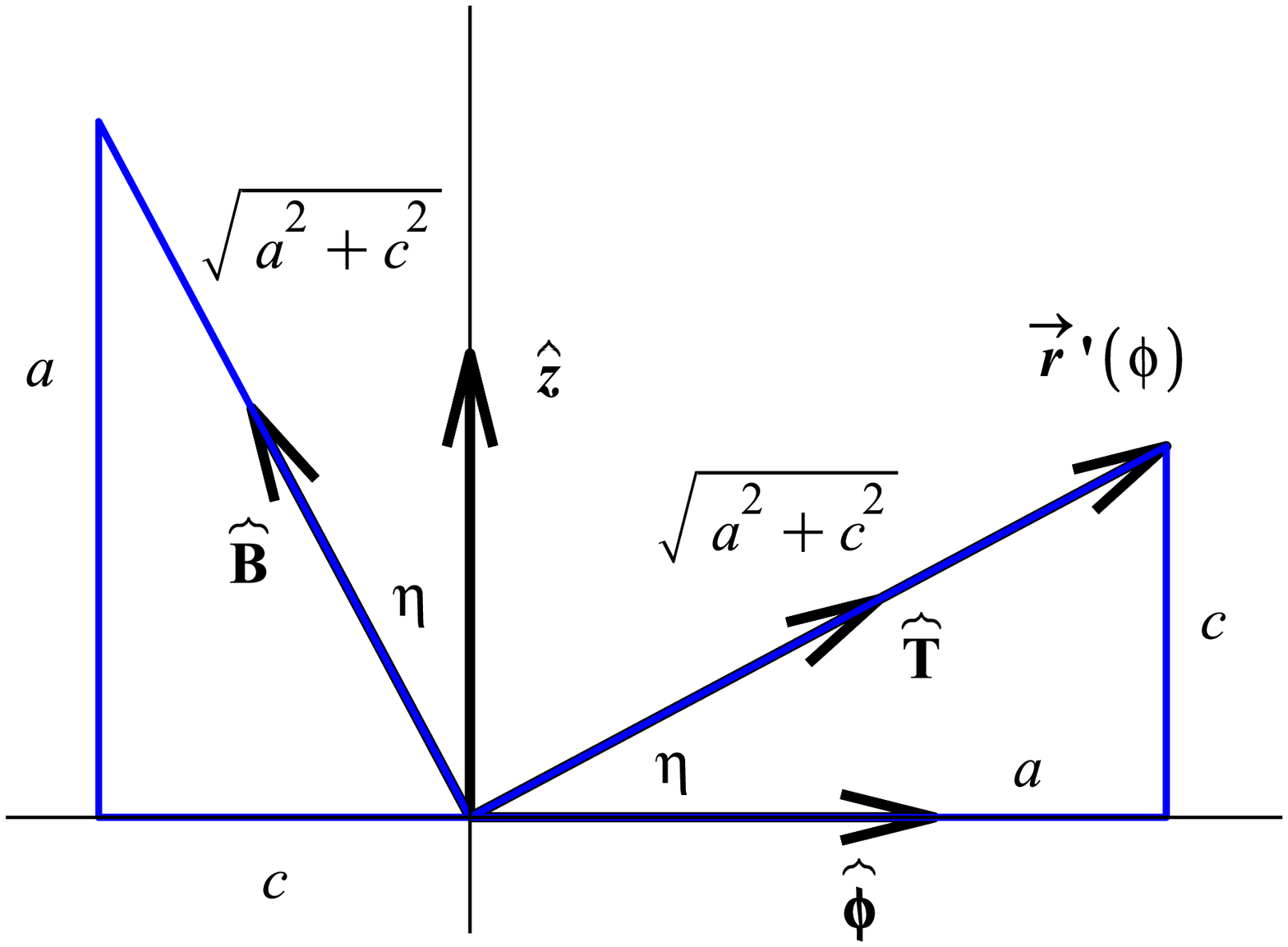}
\end{center}
\vglue0cm
\caption{
\textbf{Left:} 
A vertical half-plane cross-sectional circle of the helical tubular surface built around a helix through the center of this circle whose axis of symmetry is the $z$-axis. 
The plane of this circle is then tilted about its horizontal diameter, backwards by the angle $\psi$ from the vertical in the $\hat{\boldsymbol{\phi}}$-$\hat{\boldsymbol{z}}$ plane, and then simultaneously rotated around this axis while being translated upwards along that axis ($c>0$), so that the right hand rule wrapping fingers around the helix in the direction in which it is rising puts the thumb up. (One can also consider a left-handed helix with $c<0$.) One can introduce a radial arclength coordinate $r=b\chi=bv$ as in our previous discussions \cite{bobtorus,bobcavatappo10} but here we just use the angular parameter $\chi$ named simply $v$, while the azimuthal angle $\phi$ is designated as $u$.
\textbf{Right:} We show the tilt geometry of the special case of the orthogonally tilted cavatappo surface (``2.0") in which the tiltback angle $\psi$ equals the inclination angle $\eta$  from the horizontal of the central helix, shown in the previous diagram.
The inclination angle of the helix illustrated here is taken to be that of one version of the smooth cavatappi pasta shape: $\eta=\arctan(8/15)\approx 21.80^\circ$ \cite{blog}.
The outer ($\chi=0$, red) and inner ($\chi=\pm\pi$, green) equators are shown as bullets in the first diagram, together with the northern and southern polar helices  ($\chi=\pm\pi/2$, fuchsia).
} 
\label{fig:cross-section}
\end{figure}

After doing most of this investigation some further Googling turned up some interesting hits in computer aided design (CAD) materials written by computer people who are usually not so interested in analytical descriptions as in efficient numerical coding of surface design that has many real world applications, from industrial design to entertainment graphics, an area that has exploded with the increased computing power and software sophistication made available in recent decades. Thus one has a gap between mathematical discussions, especially those found in differential geometry texts, and newer CAD related studies. The general cavatappo 2.0 family turns out to be an example of a swept surface \cite{bishop}, namely a surface generated by moving a profile curve along a trajectory curve by specifying the orientation and possible scaling of the swept curve along its path. Freezing the scaling freedom using a circular profile curve of fixed radius leads to tube surfaces, which are best described in terms of the Frenet-Serret orthonormal frame along the trajectory curve. The orthogonally tilted cavatappo surface defined below is such a tube surface based on a helical trajectory curve, and as such has simpler intrinsic and extrinsic curvature properties.

\section{The cavatappo surface family}

The screw-symmetry about the vertical $z$-axis of ordinary Euclidean 3-space with Cartesian coordinates $(x,y,z)$ is best suited to a description starting with cylindrical coordinates
\beq\label{eq:surface}
   \begin{pmatrix} x\\ y\\ z \end{pmatrix}
  =  
  \begin{pmatrix}\rho \cos(\phi)\\\rho \sin(\phi)\\ z  \end{pmatrix}
=
  \begin{pmatrix}r \cos(\theta)\\ r\sin(\theta)\\ z  \end{pmatrix}\,,
\eeq
usually called $(\rho,\phi,z)$ by physicists and $(r,\theta,z)$ in most calculus textbooks. Let $\hat{\boldsymbol{\rho}}$, $\hat{\boldsymbol{\phi}}$ and $\hat{\mathbf{z}}$ be the unit vectors along the coordinate lines of this system, explicitly given in the following 4-parameter parametrization of a screw-symmetric surface, enlarging the previous cavatappo family by the tilt back angle $\psi$ of the plane of the circular profile curve from the vertical. This family is obtained by allowing that circular profile curve to flow along the orbits of the 1-parameter
 screw-symmetry group acting on Euclidean 3-space, corresponding to a rotation about the vertical axis accompanied by a translation along the axis by a distance proportional to the angle of rotation.
As such it is generated by a Killing vector which is a constant linear combination of the generators of translation along the $z$-axis and rotation about that axis, leading to a corresponding screw-angular momentum which is a linear combination of linear momentum along the vertical axis and angular momentum around that axis.

The parametrized ``cavatappo 2.0" surface is
\begin{eqnarray}
  \begin{pmatrix} x\\ y\\ z \end{pmatrix}
  &=&  (a+b\cos(v))\, \hat{\boldsymbol{\rho}}
     +    b\sin(v) \left( -\sin(\psi)\, \hat{\boldsymbol{\phi}} +\cos(\psi) \,\hat{\mathbf{z}}\right)
     + c u\, \hat{\mathbf{z}}
\nonumber\\
  &=&  (a+b\cos(v)) 
  \begin{pmatrix} \cos(u)\\ \sin(u)\\ 0  \end{pmatrix}
\nonumber\\
  &&\qquad  + b\sin(v)   \left( -\sin(\psi) \begin{pmatrix} -\sin(u)\\ \cos(u)\\ 0  \end{pmatrix}
                      +\cos(\psi) \begin{pmatrix} 0\\ 0\\ 1  \end{pmatrix} \right)
 + c u \begin{pmatrix} 0\\ 0\\ 1  \end{pmatrix}
\nonumber\\
 &=&   \begin{pmatrix} \left(a+b\cos(v)\right) \cos(u) +b\sin(\psi)\sin(v)\sin(u)\\
                        \left(a+b\cos(v)\right) \sin(u) -b\sin(\psi)\sin(v)\cos(u)\\
                        \phantom{\left(a+b\cos(v)\right) \sin(u)+} \kern-3ptb\cos(\psi)\sin(v)
                         + c u  \end{pmatrix}
\end{eqnarray}
where
\beq
 -\infty < u < \infty\,,\qquad  0\le v < 2\pi \mbox{\ or\ } -\pi <v \le \pi\,,
\eeq
makes these two parameters serve as coordinates on the surface (in a 1-1 relationship with points on the surface). Note that $u=\phi$ is just the usual cylindrical azimuthal coordinate, called $\theta$ in our previous investigation, and $r=b\, v$ is an arclength coordinate around the circular profile curve, which will also be referred to as the radial direction.
Adopting the same terminology common for surfaces of revolution, we call the helical orbits of the screw-symmetry group ``parallels" (the $u$ coordinate lines, the ``azimuthal direction") and the cross-sectional circles ``meridians" (the $v$ coordinate lines, the ``radial direction"). The unbounded coordinate $u$ is the azimuthal angle around the $z$-axis of symmetry, while the bounded coordinate $v$ is the angle around the cross-sectional circle.

This surface represents a cross-sectional circle (the prime meridian $u=0$) of radius $b$ with center a distance $a$ from the vertical screw rotation axis in a plane tilted by the angle $\psi$ backwards from the vertical with respect to the forward azimuthal direction (increasing $u$) which is allowed to flow along the screw symmetry orbits, which are helices. The central path of the orbit of the circle is a helix of radius $a$ and inclination angle $\eta=\arctan(c/a)$ with respect to the horizontal. 
The Killing vector generating the screw-symmetry of the Euclidean metric on $\mathbb{R}^3$ is
$\xi=a(x\,\partial_y-y\,\partial_x) + c \,\partial_z=a\,\partial_\phi + c\, \partial_z$, equal to $\partial_u$ when restricted to the surface.

Within this family is the subfamily having orthogonally tilted meridian circles where the tilt back angle $\psi$ equals the horizontal inclination angle $\eta$
\begin{eqnarray}
  \psi =\eta= \arctan\left(\frac{c}{a}\right)&&
  \quad\leftrightarrow\quad
    a\cos(\psi)=c\sin(\psi) 
\,,\nonumber \\
&&\quad\leftrightarrow\quad
   \cos(\psi)=\frac{a}{\sqrt{a^2+c^2}}\,,\ \sin(\psi) =\frac{c}{\sqrt{a^2+c^2}}
\,.
\end{eqnarray}
With this value of the inclination angle we get 
the 3-parameter orthogonally tilted cavatappo 2.0 family 
\beq
  \begin{pmatrix} x\\ y\\ z \end{pmatrix}
  =  
 \begin{pmatrix}\displaystyle \left(a+b\cos(v)\right) \cos(u) +\frac{bc}{\sqrt{a^2+c^2}}\sin(u)\\ \displaystyle 
                        \left(a+b\cos(v)\right) \sin(u) -\frac{bc}{\sqrt{a^2+c^2}}\cos(u)\\ \displaystyle 
                        \phantom{\left(a+b\cos(v)\right) b\sin(u)+\ +} 
                         \frac{ab}{\sqrt{a^2+c^2}}\sin(v) +c u  \end{pmatrix}
\,.
\eeq

Setting $\psi=0$ reduces the general cavatappo 2.0 surface to the cavatappo 1.0 surface exhaustively studied previously \cite{bobcavatappo10}, while in addition, setting $c=0$ reduces this to the torus.
George Legendre assumes the particular parameter values (we choose units of 2 mm to make a unit meridian circle $b=1$ to set the scale)
\beq
   a=\frac{3}{2}\,,\
   b=1\,,\
   c=\frac{5}{2\pi} \approx 0.7958 \sim 0.80 =\frac{4}{5}\,,\
   \eta=0\,,
\eeq 
relying on the same construction as rotational symmetry, namely starting with a plane curve in a vertical plane containing the symmetry axis to then allow to flow along the orbits of the symmetry group. However,  a little cerebral reflection suggests that restricting the model to vertical cross-sectional circles (cavatappo 1.0) is unreasonable, so adding in the tilt back angle $\psi$ fixes this limitation. Of all these tilt back choices, the orthogonal tilt back choice $\psi=\eta$ is the clear winner, rotating the cross-sectional circle so that its center moves perpendicularly to the plane containing it to sweep out all of the meridian cross-sections. Second runners up would be $\eta=0$ or $\pi/2$ so that the meridians are either vertical or horizontal. Only the sub-family of orthogonally tilted surfaces in this larger family have simpler metric properties, namely that the inner product of the tangents to the meridians and parallels is a constant.

The cavatappo 2.0 inclination angle associated with the above parameter values but using instead the very nearby rational choice $c=4/5$ 
(used by Sander Huisman in his Mathematica code for cellantani but revealed by Sol Lederman in his blog \cite{blog,NYTimes}) is
\begin{eqnarray}
  \eta &=&\arctan\left(\frac{c}{a}\right) = \arctan\left(\frac{4/5}{3/2}\right) = \arctan\left(\frac{8}{15}\right) 
\,,\nonumber\\
&& (\cos(\eta),\sin(\eta)) = \left( \frac{15}{17},\frac{8}{17} \right)
\,,
\end{eqnarray}
which is the smaller of the two angles of the right triangle whose sides and hypotenuse are serendipitously the fourth Pythagorean triplet $(8,15,17)$, leading to our canonical pasta choice
\beq
  \begin{pmatrix} x\\ y\\ z \end{pmatrix}
  =  
 \begin{pmatrix} \left(\frac{3}{2}+\cos(v)\right) \cos(u) +\frac{8}{17}\sin(u)\\[3pt]
                        \left(\frac{3}{2}+\cos(v)\right) \sin(u) -\frac{8}{17}\cos(u)\\[3pt]
                         \frac{15}{17}\sin(v)+\frac{4}{5} u  \end{pmatrix} \,.
\eeq

Of course one can easily slightly modify the profile curve as previously explored for the cavatappo rigato by adding ridges, or even rescale the circle into an ellipse in its vertical direction before tilting.
The larger family of cavatappo 2.0 surfaces, with a slight elliptical modification (replace $b\,\sin(v)$ in Eq.~[\ref{eq:surface}] with $d\,\sin(v)$ to get an elliptical cross-section), allows another less esoteric application, of pedagogical use for students learning about special relativity. By re-interpreting the vertical symmetry axis as the time axis in a 3-dimensional flat Minkowski spacetime in inertial coordinates, one can adjust the tilt back angle to rotate the cross-sectional plane to a horizontal orientation and Lorentz contract one axis of the circular cross-sections so that they appear to be circular in the rest frame of the central world line in order that the sphere undergoes Born rigid motion, thus appearing Lorentz contracted in the inertial frame in which the axis is at rest \cite{bobfermi}. However, in this case it is much easier to simply redo the orthogonality directly in the flat Lorentzian geometry of 3-dimensional Minkowski spacetime, simply boosting a circular horizontal cross-section to the local rest frame of a timelike helical world line.

To describe this surface as a swept surface, one needs the Serret-Frenet frame along the central helix consisting of the unit tangent $\hat T$, the unit normal $\hat N$ and the unit binormal $\hat B$
\begin{eqnarray}
\vec r (\phi) &=& \langle a\cos\phi,a\sin\phi,c \phi \rangle = a \hat\rho + c \phi \hat z \,,\nonumber\\
\vec r\, '(\phi) &=& \langle -a\sin\phi,a\cos\phi,c \phi \rangle = a \hat\phi + c \hat z \,,\nonumber\\
\hat T (\phi) &=& (a^2+c^2)^{-1/2} (a\hat\phi +c\hat z) \,,\nonumber\\
\hat N (\phi) &=& -\hat \rho \,,\nonumber\\
\hat B (\phi) &=& (a^2+c^2)^{-1/2} (-c\hat\phi +a\hat z) \,.
\end{eqnarray}
One can introduce an arclength parametrization $s=(a^2+c^2)^{1/2} \phi$, and easily evaluate the radius of curvature 
$\mathcal{R}=\kappa^{-1} = a(1+c^2/a^2)$.
The orthogonally tilted cavatappo 2.0 surface is a tube surface consists of a circle of radius $b$ about the origin in each normal plane spanned by $\hat N$ and $\hat B$, typically parametrized by the usual polar angle $\Theta$ in that plane with ordered basis $\{\hat N,\hat B\}$.
\beq
\vec r (\phi,\Theta) = \vec r (\phi)+ b\left(\cos\Theta\, \hat N(\phi) + \sin\Theta \, \hat B (\phi)\right)
\eeq
where $\Theta=\pi-\chi$,
while the general surface is a swept surface with an additional rotation of the circular profile curve by the angle $\eta-\psi$ about the normal vector (right hand rule orientation).

\section{The metric}

Taking the differentials of the general surface parametrization of $\langle x,y,z\rangle$, inserting them into the Euclidean metric
$ds^2=dx^2+dy^2+ dz^2$
%\epsilon dz^2$, where $\epsilon=\pm 1$ allows one to consider the vertical axis as an inertial time in 3-dimensional Minkowski spacetime, 
and manipulating the result leads to 
\begin{eqnarray}\label{eq:metric}
ds^2
&=&\left( (a+b\cos(v))^2 +b^2 \sin^2(v)\sin^2(\eta)+c^2\right)\, du^2\nonumber\\
  &&\qquad    + 2b\left( \strut\left(-a\sin(\eta)+c\cos(\eta)\right)\cos(v)  -b\sin(\eta)   \right) \, du \, dv\nonumber\\
  &&\qquad    + b^2 dv^2
\,.
\end{eqnarray}
For the orthogonally tilted cavatappo 2.0 subfamily this reduces to
\begin{eqnarray}
%ds^2&=&\left(\strut (a+b\cos(v))^2 +\frac{b^2c^2}{a^2+c^2} \sin^2(v) +c^2\right)\, du^2
%\nonumber\\
%  &&\qquad    - \frac{2b^2c}{\sqrt{a^2+c^2}}  \, du \, dv + b^2 dv^2
%\nonumber\\
ds^2&=&\frac{ (a^2+c^2 +ab\cos(v))^2 +b^2 c^2}{a^2+c^2}du^2\nonumber\\
  &&\qquad    - \frac{2b^2c}{\sqrt{a^2+c^2}}  \, du \, dv + b^2 dv^2
\,,
\end{eqnarray}
and with the canonical parameter values
\begin{eqnarray}
ds^2&=&%\left( \left(\frac{3}{2}+\cos(v)\right)^2 +\left(\frac{8}{17}\right)^2 \sin^2(v) +\left(\frac{4}{5}\right)^2\right)\, du^2\nonumber\\
  \frac{ \left( \left(\frac{3}{2}\right)^2 + \left(\frac{5}{4}\right)^2  + \left(\frac{3}{2}\right) \cos(v)\right)^2}{
         \left(\frac{3}{2}\right)^2 + \left(\frac{5}{4}\right)^2 } \, du^2
%\nonumber\\   &&\qquad   
 - \frac{16}{17}  \, du \, dv + dv^2\,.
\end{eqnarray}
Note that $r=b\, v$ is an arclength coordinate around the meridians measured from the outer equator, to be referred to as the radial direction as in the surface of revolution case. The azimuthal angle $u$ coincides with the cylindrical coordinate usually referred to as $\phi$ by physicists, but the variable names $(r,\theta)=(u,v)$ correspond to the previous discussion of the torus and cavatappi 1.0 family.

One can complete the square on either of the squares in these quadratic forms to achieve two different orthogonal canonical forms for the metric analogous to adapting the metric to the orbits of the time translation symmetry group in a stationary spacetime \cite{mfg}. For simplicity we only give the results for the cavatappo 2.0 subfamily and the special values given above. For example completing the square on $dv$  yields in general and for the canonical parameter values
\begin{eqnarray}
ds^2 &=& \left( (a+b\cos(v))^2-\frac{b^2 c^2 \cos^2(v)}{a^2+c^2} +c^2\right) \, du^2
        + b^2 \left(dv-\frac{c\, du}{\sqrt{a^2+c^2}}  \right)^2
\nonumber\\
 &=& \left( \left(\frac32+\cos(v)\right)^2-\left(\frac{8}{17}\right)^2 \cos^2(v) +\left(\frac{4}{5}\right)^2 \right) \, du^2\,,
\nonumber\\
&&\qquad
        +  \left(dv-\frac{8}{17} du  \right)^2\,.
\end{eqnarray}

It is helpful to introduce some orthogonal decomposition quantities for both cases, with appropriate notation for each.
Adapting the decomposition to be orthogonal with respect to the symmetry group orbits, namely the $u$ coordinate lines or  ``threads," gives the threading decomposition (complete the square on $du$)
\begin{eqnarray} \label{eq:metricM}
 ds^2
&=& g_{uu} \, du^2 +2 g_{uv} \, du\, dv + g_{vv}\, dv^2\\
&=& M^2 (du+M_v dv)^2 +\gamma_{vv} \, dv^2
\,,\nonumber\\
&=& (\omega^{\top})^2 +(\omega^{\hat v})^2 \,,
\end{eqnarray}
which requires computation of the following quantities:
\begin{eqnarray}
  M &=&(g_{uu})^{1/2}
= \left(\strut (a+b\cos(v))^2 +\frac{b^2c^2}{a^2+c^2} \sin^2(v) +c^2\right)^{1/2}
\nonumber\\
&=&  \left(\frac{\strut (a^2+c^2+ab\cos(v))^2 +b^2c^2}{a^2+c^2}\right)^{1/2}
\nonumber\\
&=&  \frac{10}{17} \left( \left(\frac32\cos(v) +\left(\frac{17}{10}\right)^2 \right)^2 +\left(\frac{4}{5}\right)^2 \right)^{1/2}
%&=& \left( \left(\frac32+\cos(v)\right)^2+\left(\frac{8}{17}\right)^2 \sin^2(v) +\left(\frac{4}{5}\right)^2 \right)^{1/2}
\,,
\nonumber\\     
M_v &=& \frac{g_{uv}}{M^2} =  -\frac{b^2c}{M^2\sqrt{a^2+c^2}} 
     = -\frac{8}{17 M^2} 
\nonumber\\
   M^v &=& M_v/g_{vv}= -\frac{c}{M^2\sqrt{a^2+c^2}}=-\frac{8}{17 M^2} 
\,,
\nonumber\\
  \gamma_{vv} &=& g_{vv} -g_{uv}^2/M^2 
  =
%\frac{b^2 (a^2+c^2+ab\cos(v))^2}{\left((a^2+c^2)(a+b\cos(v))^2+c^2+b^2c^2\sin^2(v) \right)}
\frac{b^2 (a^2+c^2+ab\cos(v))^2}{\strut (a^2+c^2+ab\cos(v))^2 +b^2c^2}
\nonumber\\
&=& \frac{  \left(\frac{17}{10}\right)^2 \left(\left(\frac{17}{10}\right)^2 +\frac{3}{2}\cos(v) \right)^2}
%     {\left(\frac{3}{2}+\cos(v)\right)^2 +\left(\frac{8}{17}\right)^2 \sin^2(v)+\left(\frac{4}{5}\right)^2}
     {\left(\frac{17}{10}\right)^2 \left(\left(\frac{17}{10}\right)^2 +\frac{3}{2}\cos(v) \right)^2 +\left(\frac{4}{5}\right)^2 }
\,.
\end{eqnarray}
This corresponds to the orthonormal frame and dual frame
\begin{eqnarray}
&& \Vec{e}_{\top} = M^{-1/2} \partial_u\,,\
   \VecS{\epsilon}_{\hat v} = \gamma_{vv}^{-1/2} (\partial_v-M^v\partial_u)\,,
\nonumber\\   
&& \omega^{\top} = M (du+M_v\, dv)\,,\
   \omega^{\hat v} = \gamma_{vv}^{1/2} \omega^v\,.
\end{eqnarray}
This decomposition is useful for dealing with the screw-angular momentum which is a constant quantity along geodesics.

The alternative is the slicing decomposition which is orthogonal with respect to the ``slices" of constant $u$, namely the $v$ coordinate lines or meridians (complete the square on $dv$)
\begin{eqnarray}\label{eq:metricN}
 ds^2 &=& N^2 du^2 +g_{vv} \, (dv+N^v du)^2
\,,\nonumber\\
&=& (\theta^\bot)^2 +(\theta^{\hat v})^2 \,,
\end{eqnarray}
which requires computation of the following quantities:
\begin{eqnarray}
  N &=&(g_{uu}-g_{uv}^2/g_{vv})^{1/2} 
%= \left( (a+b\cos(v))^2-\frac{b^2 c^2 \cos^2(v)}{a^2+c^2} +c^2\right)^{1/2} 
=  \frac{a^2+c^2 +a b \cos(v)}{\sqrt{a^2+c^2}} 
\nonumber\\
 &=& \left(\frac32+\cos(v)\right)^2-\left(\frac{8}{17}\right)^2 \cos^2(v) +\left(\frac{4}{5}\right)
\,,\nonumber\\ 
    N^v &=& \frac{g_{uv}}{g_{vv}} =   -\frac{c}{\sqrt{a^2+c^2}} = -\frac{8}{17}
\,,\ 
\nonumber\\
N_v &=& g_{uv}=g_{vv} N^v=-\frac{b^2 c}{\sqrt{a^2+c^2}} =-\frac{8}{17}
\,,
\nonumber\\
g_{vv} &=& b^2 =1
\,.
\end{eqnarray}
This corresponds to the orthonormal frame and dual frame
\begin{eqnarray}
&& e_{\bot} = N^{-1/2} (\partial_u-N^v\partial_v)\,,\
   e_{\hat v} = g_{vv}^{-1/2} \partial_v\,,
\nonumber\\   
&& \omega^{\top} = N \, du \,,\
   \omega^{\hat v} = g_{vv}^{1/2} (dv+N^v du)\,.
\end{eqnarray}
This latter choice is useful for the cavatappo 2.0 family where the meridians admit a new simple threading orthogonal to them, yielding an orthogonal grid on the surface discussed below. 
 
The square root of the metric determinant determines the surface element $dS= \det(g)^{1/2}du\,dv$, and is easily seen from these orthogonal decompositions as the square root of the product of the diagonalized metric coefficients
\begin{eqnarray}
  \det(g)^{1/2} &=& N g_{vv}^{1/2} = M \gamma_{vv}^{1/2} 
    %=  b \left( (a+b\cos(v))^2-\frac{b^2 c^2 \cos^2(v)}{a^2+c^2} +c^2\right)^{1/2} 
\nonumber\\
 &=& \frac{b}{\sqrt{a^2+c^2}} (a^2+c^2 +a b \cos(v))
  = \frac{17}{10} + \frac{15}{17} \cos(v)\,.
\end{eqnarray}
This is easily integrated for one revolution of the surface to yield
$$
\int_0^{2\pi} \int_0^{2\pi}  \det(g)^{1/2}\, dv\,du = (2\pi b)(2\pi\sqrt{a^2+c^2})
$$
which is the product of the circumference $C=2\pi b$ of the orthogonal circular cross-sections and the arclength $L=2\pi \sqrt{a^2+c^2}$  of the central curve for one revolution of the symmetry axis, an example of the theorem of Pappas for surfaces of revolution like the torus extended to the screw-symmetric case. Note that in the absence of the orthogonal tilt condition, this integral results in a many screen long formula dense in elliptic functions, showing that the orthogonality condition results in a much simpler geometry than the general case.

Note that it is not easy to find the compact forms of these various expressions in Eqs.~[\ref{eq:metric}], [\ref{eq:metricM}, [\ref{eq:metricN}] for $g_{uu}$, $N$, $M$ and $\gamma_{vv}$ since computer algebra systems do not have such simplication tools.

\begin{figure}
\typeout{*** EPS figure cavatappo threading splitting}
\begin{center}
%\vglue-2cm
{\includegraphics[scale=0.4]{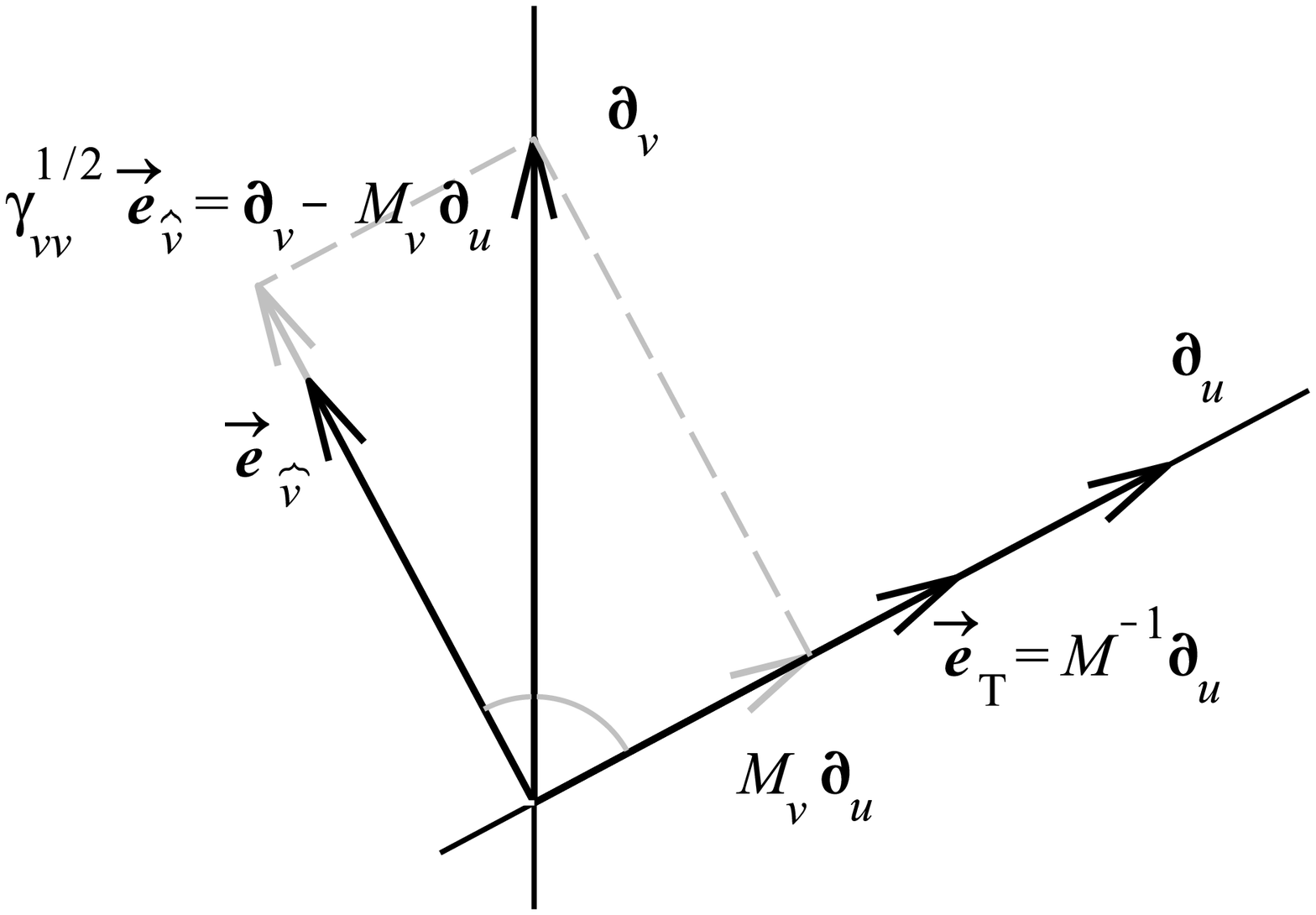}} \vglue-1cm
{\includegraphics[scale=0.4]{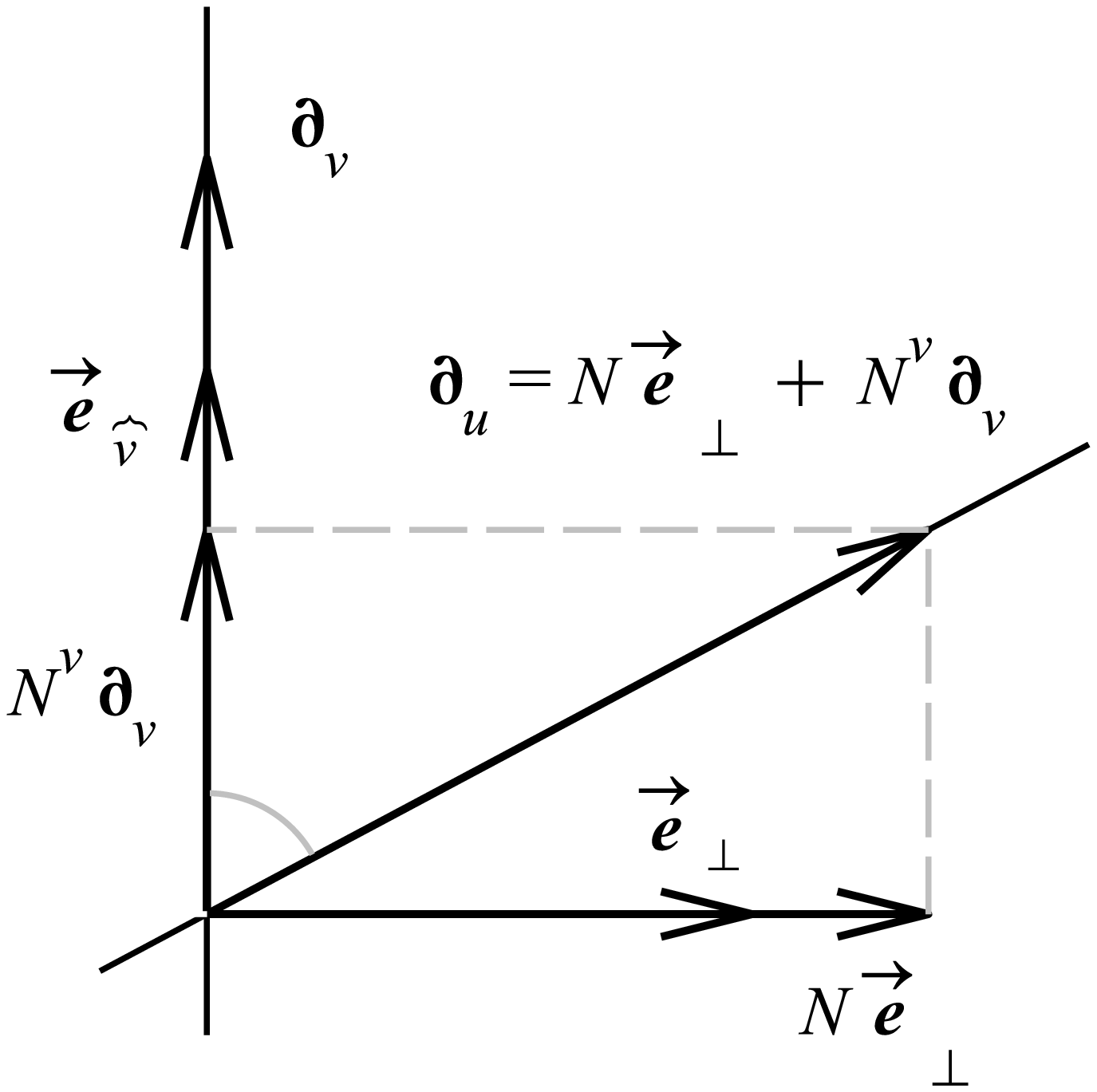}}
\end{center}
\vglue-1cm
%\vglue-2cm

\caption{
\textbf{Above:}
The orthogonal ``threading" decomposition of the tangent space with respect to the $u$ coordinate lines (parallels) in the Euclidean case.
\textbf{Below:}
The corresponding ``slicing" decomposition with respect to the $v$ coordinate lines (meridians).
}\label{fig:threadingsplit}
\end{figure}

The scalar curvature (twice the Gaussian curvature) is easily calculated with a computer algebra system
\beq
    R=2K= \frac{2a \cos(v)}{b(a^2+c^2+ab\cos(v))}
\eeq
with extreme values at the equators $v=0,\pi$
\beq
   R_{\rm ext}  = 2K_{\rm ext}=\pm \frac{2a}{b(a^2+c^2\pm ab)}\,.
\eeq
$R$ vanishes at the two polar helices, which separate the regions of positive and negative curvature. Since $K^{-1}= (k_1k_2)^{-1} = r_1
r_2$ has the interpretation as the product of the inverses of the two principal curvatures which are the radius of curvature of the orthogonal curves along the principle directions of the extrinsic curvature (second fundamental form), and $b$ is the radius of curvature of the meridians, then
\beq
  \mathcal{R} =
   \frac{a^2+c^2+ab\cos(v)}{a \cos(v)}
\eeq
must be the radius of curvature of the orthogonal trajectories to the meridians pictured in Fig.~\ref{fig:grids}, which are lines of curvature. 
In fact when $c=0$ this reduces to 
\beq
    \mathcal{R} =
   \frac{a+b\cos(v)}{ \cos(v)} =  \frac{R(v)}{ \cos(v)}\,,
\eeq
in terms of the radius $R(v)$ of the parallels of the torus, which are horizontal azimuthal circles, projected onto the normal plane along the parallels. This is verified by an explicit calculation of the extrinsic curvature tensor and its eigenvalues.

\section{Geodesic equations}

The second order geodesic equations for affinely parametrized geodesics $u(\lambda)$, $v(\lambda)$ are found to be
\def\diff(#1, #2){\frac{d #1}{d #2}}
\begin{eqnarray}
&&
\mathcal{D}  \frac {d^{2}u}{d{\lambda }^{2}}
+\frac{a b c\sin(v)}{\sqrt{a^2+c^2}} \left(\diff(u, t)\right)^2 - 2  a b \sin(v) \diff(u, t) \diff(v, t) = 0\,,
\nonumber\\
&&
\mathcal{D} {\frac {d^{2}v}{d{\lambda }^{2}}}
+ \frac{ a \sin(v) \left((a^2+c^2+ab\cos(v))^2+b^2c^2\right)}{b (a^2+c^2) } \left(\diff(u, \lambda)\right)^2
\nonumber\\
&&\qquad
-\frac{2 a b c \sin(v)}{\sqrt{a^2+c^2}} \diff(u, \lambda) \diff(v, \lambda)= 0
\nonumber\\
&&\qquad\mbox{where\quad} \mathcal{D} = a^2+c^2+a b \cos(v)\,.
\end{eqnarray}
One obvious consequence of these equations is that $du/d\lambda =0$ reduces these equations to $d^2v/d\lambda^2=0$ making the radial coordinate $v$ along the meridians a linear function of the affine parameter $\lambda$ for constant values of the azimuthal coordinate $u$, so that the meridian circles are geodesics in this new family. Similarly $\sin(v)=0$ which implies $v=0,\pi$ shows that the inner and outer equators are also geodesics.

\begin{figure}
\typeout{*** EPS figure cavatappo}
\begin{center}
\vglue-2cm
\includegraphics[scale=0.40]{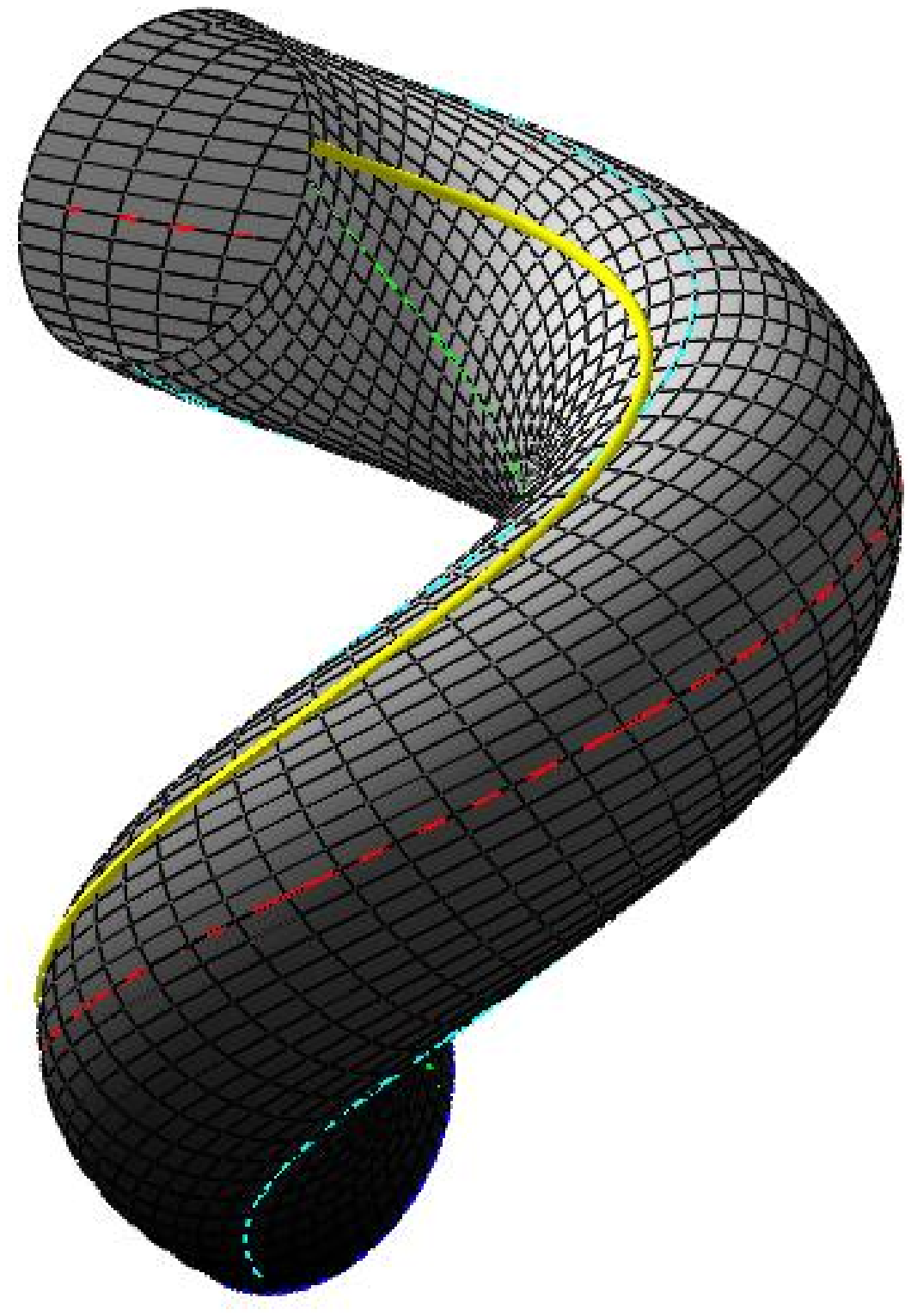}\kern-1cm 
\includegraphics[scale=0.40]{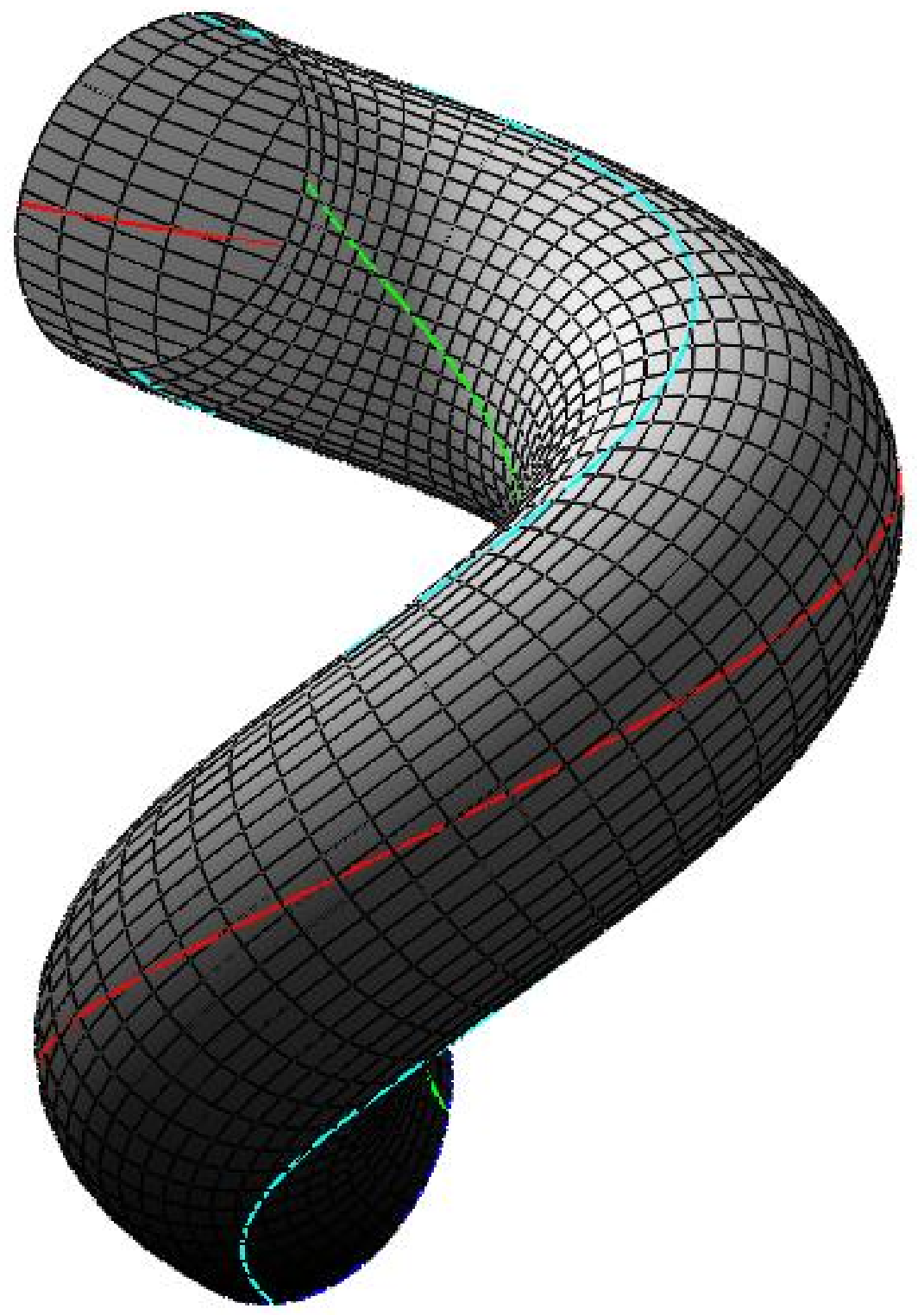}
\end{center}
\vglue-2cm

\caption{The original (left) and new (right) orthogonal grids on the cavatappo 2.0 surface, shown with a red outer equator and green inner equator and cyan polar circles. Both grids have the same meridians, but the parallels in the original grid are replaced by orthogonal curves (to the meridians) which spiral around the tubular surface. To show the comparison one spiral is shown in yellow in the original grid.
}\label{fig:grids}
\end{figure}

The orthogonality condition for the new parallels which are the orthogonal trajectories to the meridians is
\beq 
 g_{vv}^{-1/2}\omega^{\hat v} =  dv+N^v du 
  = dv -\frac{b^2 c}{\sqrt{a^2+c^2}} du =d\left(v-\frac{b^2 c}{\sqrt{a^2+c^2}} u\right) =0
\eeq
so that the following new radial coordinate is constant along the new orthogonal parallels, both in general and for our canonical parameter values
\beq
  v_\bot \equiv v-\frac{b^2 c}{\sqrt{a^2+c^2}} u =  v-\frac{8}{17} u
\,.
\eeq
Substituting $v=v_\bot+(8/17) u$ into the original surface parametrization yields an orthogonal coordinate grid on the surface consisting of meridians and helices which spiral around the tubular surface in the positive radial direction.
Fig.~\ref{fig:grids} compares the two grids. Neither set of parallels are geodesics. Evaluating the extrinsic curvature or shape tensor for the surface confirms that this new orthogonal grid consists of the lines of curvature and indeed the two principle curvatures are $1/b$ for the meridians, and $1/\mathcal{R}$ for their orthogonal trajectories.

However, the curves orthogonal to the parallels, characterized by zero screw-angular momentum, are geodesics.
\begin{eqnarray}
  && 0 =  M^{-1}\omega^{\top} = du-M_v\, dv  \rightarrow
\nonumber\\
  &&u=\int_0^v -M_v dv
   = \int_0^v  \frac{8}{17 M^2} dv\,.
\end{eqnarray}
This is exactly integrable but fills many screens with the result. For the canonical values, the increment in azimuthal angle during one radial revolution is
\beq
\Delta u = \int_0^{2\pi} -M_v dv \approx 1.3738\approx 78^o \approx 0.2186\,\mbox{revs}\,.
\eeq
Thus there are about 5 wrappings of the tubular surface during 1 revolution about the vertical axis.

\section{Constants of the motion}

If $(u(\lambda),v(\lambda))$ is an affinely parametrized geodesic of this metric on the surface, its tangent
\beq
 U =   \frac{du}{d\lambda} \partial_u + \frac{dv}{d\lambda} \partial_v
  = U^{\hat u} e_{\top} + U^{\hat v} \epsilon_{\hat v}
\eeq
can be expressed in terms of its components in either of the orthogonal surface frames. In particular in the threading decomposition one has
\begin{eqnarray}
  U^{\hat u} = M \left(  \frac{du}{d\lambda} + M_v \frac{dv}{d\lambda} \right ) 
\quad %\nonumber\\
  U^{\hat v} = \gamma_{vv}^{1/2}  \frac{dv}{d\lambda} \,,
\end{eqnarray}
This can be inverted to yield
\begin{eqnarray}  \label{eq:velocityrelation}
 U^u = \frac{du}{d\lambda} = M^{-1} U^{\hat u} - M_v \gamma_{vv}^{-1/2}   U^{\hat v}\,,
\quad %\nonumber\\
  U^v= \frac{dv}{d\lambda} = \gamma_{vv}^{-1/2}  U^{\hat v}
\,. 
\end{eqnarray}

\begin{figure}[p]
\typeout{*** EPS figure cavatappo}
\begin{center}
%\vglue-2cm
 \includegraphics[scale=0.30]{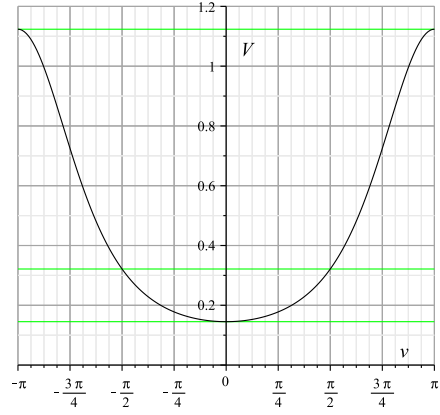}
\end{center}
%\vglue-2cm

\caption{The screw-centrifugal potential. Shown in green are the energy levels of the inner and outer equators and of the northern and southern polar helices.
}\label{fig:screw-centrifugalpot}
\end{figure}

The component of the tangent vector along the Killing vector field  (the conserved screw-angular momentum) is a constant along the geodesic
\beq
  \ell  =U_u = g_{uu}U^u+g_{uv}U^v
   = M U^{\hat u} = M^2 \left( \frac{du}{d\lambda}+M_v \frac{dv}{d\lambda} \right)
\,, 
\eeq
as is its length, half of which we call the energy (re-expressing the angular velocity term in terms of the conserved angular-momentum)
\beq
   {\frac12}\left( U^{\hat u} \right)^2 +   {\frac12}\left( U^{\hat v} \right)^2
 =  {\frac12}\left( U^{\hat v} \right)^2 +  {\frac12}\frac{\ell ^2}{M^2} =  E
\eeq
so
\beq
   V =  {\frac12}\frac{\ell ^2}{M^2}
=\frac{(a^2+c^2)\ell ^2}{2(a^2+c^2+ab\cos(v))^2 +b^2c^2}
\eeq
acts as an effective potential for the radial motion. If the screw-angular momentum is nonzero, then we might as well set it equal to 1 and use the energy parameter to distinguish initial data.
Fig.~\ref{fig:screw-centrifugalpot} shows the nonzero screw-angular momentum potential for the orthogonally tilted cavatappi 2.0 pasta surface for the canonical parameter values.

\begin{figure}
\typeout{*** EPS figure cavatappo timelike initial data}
\begin{center}
%\vglue-2cm
\includegraphics[scale=0.35]{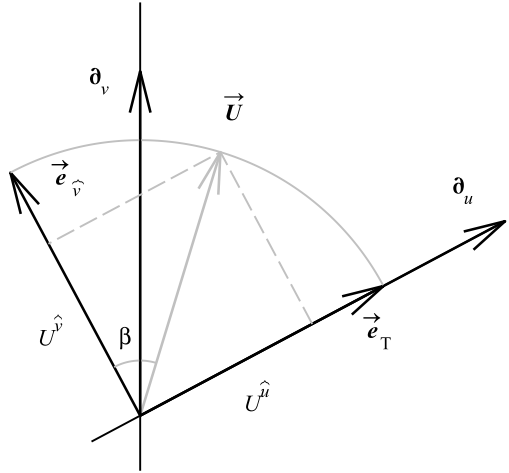}
\end{center}
%\vglue-2cm

\caption{Orthogonal decomposition of initial data in the vertical tangent plane at the outer equator for a unit tangent vector $\Vec{U}$ with respect to the orthonormal frame adapted to the parallels along $\partial_u$.
}\label{fig:initialdata}
\end{figure}

\section{Initial value problem for nonzero screw-angular momentum}

For interpretation it is useful to pose initial data for numerical solution of the geodesic equations at the origin of coordinates on the outer equator at the prime meridian in terms of an angle $\beta$ with respect to the threading orthonormal frame, namely
\beq
 \Vec{U}   = U^{\hat u} \Vec{e}_{\top} + U^{\hat v} \VecS{\epsilon}_{\hat v}\,,\quad
 (U^{\hat u} , U^{\hat v}) = 2E (\sin\beta,\cos\beta) \,.
\eeq
We choose the angle from the meridian direction of increasing $v$ to be consistent with in our previous cavatappo 1.0 discussion, which in turn was established for the special case of the torus problem.
Then to produce initial data for the second order differential equations we need the coordinate velocities
\begin{eqnarray}  \label{eq:velocityrelationinitial}
\frac{du}{d\lambda} = 2E\left( M^{-1} \sin\beta - M_v \gamma_{vv}^{-1/2}   \cos\beta\right)\,,
\quad
\frac{dv}{d\lambda} = \gamma_{vv}^{-1/2}  \cos\beta
\,, 
\end{eqnarray}
in terms of which the conserved screw-angular momentum is
\beq
  \ell = M (2E) \sin\beta \,.
\eeq
For numerical solution a unit arclength parametrization $2E=1$ is simplest, making $U$ the unit tangent.

\section{Classification of geodesics}

As on the torus, the geodesics fall into two classes: those radially unbounded orbits which wrap around the tubular surface, crossing the inner equator, and those which are radially bounded in their motion, using the physics perspective in interpreting tracing out the geodesics as motion constrained to lie on the surface.
One can also investigate the geodesics on this surface which are periodic with respect to the azimuthal coordinate, namely those for which the radial variable $v$ undergoes an integral number of radial oscillations during an integral number of azimuthal revolutions. These geodesics, when projected onto a horizontal plane, are closed curves, while the rest are not. This mirrors the approach that has already been exhaustively applied to the special case of tori \cite{bobtorus} and which can be carried over to the present discussion in parallel.

\section{From Euclidean to Lorentzian helices: the central helix}

By changing the signature of the flat metric on 3-space, we can model interesting special relativistic phenomena while providing a fascinating toy  model for general relativistic calculations.
This requires suppressing one spatial dimension of the 4-dimensional Minkowski spacetime, for example, consider the horizontal equatorial circle on a classical electron in a circular orbit in a horizontal plane at constant speed and ignore the vertical dimension ($z$). The horizontal component of its spin vector in the $x$-$y$ plane  undergoes a Thomas precession in that plane due to the twisting of its world line in spacetime \cite{mtw}. This is easily described in a 3-dimensional Minkowski space with inertial coordinates $(x,y,t)\equiv (x^1,x^2,x^0)$. If we model the electron as undergoing Born rigid motion, the only natural choice, so that it is a sphere in its local rest space, then suppressing the electron's vertical dimension, its equatorial plane cross-section in its local rest frame is a circle, which undergoes Lorentz contraction as seen in the inertial coordinates.

\begin{figure}
\typeout{*** EPS figure cavatappo triangles tilt diagram}
\begin{center}
\vglue-1cm
 \includegraphics[scale=0.35]{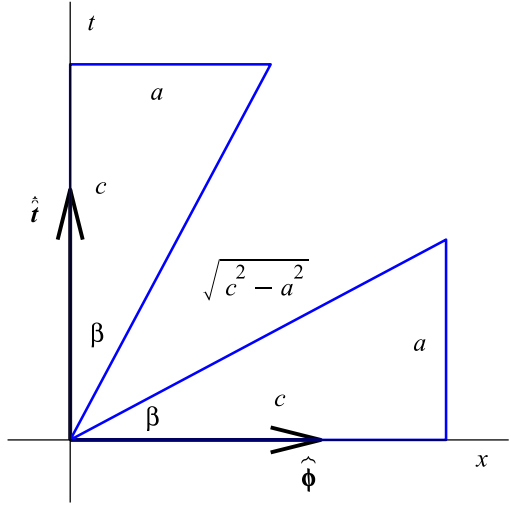}
\end{center}
%\vglue-2cm

\caption{The Lorentz tilt diagram showing the inclination of the central helix (upper triangle) and the inclination of the orthogonal local rest space (lower triangle). The timelike condition for the speed $c/a>1$ stretches out the helix along its axis, with the hyperbolic inclination angle $\beta = {\rm arctanh}(a/c)$ and Euclidean inclination angle $\eta=\arctan(a/c)$, relative to the horizontal direction.
}\label{fig:lorentztilt}
\end{figure}

We rename the vertical coordinate $z$ to $t$ so that the orthonormal Cartesian coordinates on Euclidean 3-space become inertial coordinates on Lorentzian 3-dimensional Minkowski spacetime, but keep it in the third position so that in our technology plots the time direction is vertical, and we use polar coordinates in the 2-dimensional plane spaces to describe the helical motion
\beq\label{eq:surfaceLorentz}
   \begin{pmatrix} x\\ y\\ t \end{pmatrix}
  =  
  \begin{pmatrix}\rho \cos(\phi)\\\rho \sin(\phi)\\ t  \end{pmatrix}
=
  \begin{pmatrix}r \cos(\theta)\\ r\sin(\theta)\\ t  \end{pmatrix}\,,
\eeq
with the Minkowski metric, letting $i,j=0,1,2$, with inertial coordinates (orthonormal Cartesian coordinates) $(x^1,x^2,x^0)=(x,y,t)$,
\beq
  ds^2=\eta_{ij}dx^i dx^j= dx^2+dy^2-dt^2 =  d\rho^2+\rho^2 d\phi^2 -dt^2 \,,
\eeq
and metric matrix
\beq
  (\eta_{ij}) = \begin{pmatrix} 1&0&0\\ 0&1&0\\ 0&0&-1\end{pmatrix} = (\eta^{ij})\,.
\eeq
We take a timelike helical world line $\rho = a, \phi=\omega t$ but parametrized by the proper time $\tau$ (spacetime arc length for a timelike curve), and assuming $c>0$ for simplicity of discussion
\begin{eqnarray}
  \vec{\mathbf{x}}=   \langle x,y,t\rangle
 &=&  \langle a \cos(\phi),a \sin(\phi), c\phi\rangle
 = \langle a \cos(\omega t), a \sin(\omega t), t\rangle
\nonumber\\
 &=& \langle a \cos(\gamma\omega \tau),a \sin(\gamma\omega \tau), \gamma\tau\rangle \,,
\end{eqnarray}
where the inertial time angular velocity and proper time angular velocities of the circular motion are related by
\beq
\frac{d\phi}{dt}=  \omega = \frac{1}{c}\,,\ 
 \frac{d\phi}{d\tau} = \gamma\omega =\frac{1}{\sqrt{c^2-a^2}}
\,.
\eeq
The azimuthal speed (not to be confused later with the surface parameter $v$) is
\beq
  v=\left( \left(\frac{dx}{dt}\right)^2 + \left(\frac{dy}{dt}\right)^2 \right)^{1/2} = a\frac{d\phi}{dt}
   = a\omega = \frac{a}{c} \equiv \tanh\beta\,,
\eeq
defining the rapidity $\beta$, with associated gamma factor
\beq
 \gamma=(1-v^2)^{-1/2} = \frac{c}{\sqrt{c^2-a^2}} = \cosh\beta \,, \qquad \frac{d\tau}{dt}=\gamma
\eeq
relating inertial and proper times along the world line, and
\beq
  \sinh\beta=\gamma v = \frac{a}{\sqrt{c^2-a^2}} \,.
\eeq 

This world line has a
timelike unit tangent vector (the 4-velocity)
\beq
  \vec{\mathbf{U}}= \frac{d}{d\tau} \langle x,y,t\rangle
 = \langle -a \gamma\omega \sin(\gamma\omega \tau),a \gamma \omega \cos(\gamma\omega \tau), c\rangle
\eeq
with a helical inclination slope $c/a$ related to the gamma factor and hyperbolic angle of inclination $\beta$ and the Euclidean angle $\eta$ from the horizontal by
\begin{eqnarray}
  &&\frac{a}{c} \equiv \cot\eta \equiv \tanh\beta =v 
\,.
\end{eqnarray}
The hyperbolic angle $\beta$ is the rapidity associated with the speed of the world line, and the azimuthal direction boost parameter from the inertial coordinate frame to the comoving local rest space along the world line $LRS_U$, namely the spacelike plane orthogonal to the tangent vector $\vec{\mathbf{U}}$ to the helix at each point on that helix.

The total arclength along this world line  for one revolution is the proper time elapsed along it, namely the proper period $T_o$. The inertial time period is determined by $\omega T = 2\pi$, so $T=2\pi/\omega=2\pi c$. The proper period is just related to the inertial time by time dilation, namely $T_0 = T/\gamma= 2\pi \sqrt{c^2-a^2}$ (since it is determined by the relation $\omega_o T_o = 2\pi$).

The Frenet-Serret approach to describing the geometry of a parametrized curve in Euclidean space is easily carried over to timelike world lines in Minkowski spacetime \cite{bobdgn}. First discussed for general spacetimes by Synge \cite{synge}, it later was applied to world lines in 4-dimensional Minkowski spacetime (including timelike helices as a very special case) to study trajectories of particles in homogeneous electromagnetic fields \cite{HSV} and later applied to study spin precession along circular orbits in black hole spacetimes \cite{IV}.

The unit normal $\vec{\mathbf{N}}$ and the 4-acceleration  $\vec{\mathbf{A}}$ are defined as in the 3-dimension\-al Euclidean case by
\begin{eqnarray}
  \vec{\mathbf{A}}  &=&  \frac{d}{d\tau} \vec{\mathbf{U}}= \frac{d^2}{d\tau^2} \langle x,y,t\rangle
 = \langle -a \gamma^2 \omega^2 \cos(\gamma\omega \tau),- a \gamma^2 \omega^2 \sin(\gamma\omega \tau), 0\rangle
\nonumber\\ &=&
  \underbrace{ a \gamma^2 \omega^2}_{\displaystyle \kappa }  \underbrace{ \langle - \cos(\gamma\omega \tau),-  \sin(\gamma\omega \tau), 0\rangle}_{\displaystyle\vec{\mathbf{N}} } \,,
\end{eqnarray}
with curvature $\kappa\ge0$ equal to the magnitude of the acceleration.
The binormal $\vec{\mathbf{B}}$ is the direction of the component of the derivative of the unit normal orthogonal to the unit tangent, completing the Lorentz signature Frenet-Serret relations for the timelike world line
\beq
  \frac{d}{d\tau} \vec{\mathbf{U}} = \kappa \vec{\mathbf{N}}\,,\
  \frac{d}{d\tau} \vec{\mathbf{N}} = \kappa  \vec{\mathbf{U}} +\Omega  \vec{\mathbf{B}}\,,\
  \frac{d}{d\tau} \vec{\mathbf{B}} = -\Omega  \vec{\mathbf{N}}\,.
\eeq
The two spacelike normal vectors $\vec{\mathbf{N}}$ and $\vec{\mathbf{B}}$ span the local rest space.
We can get the binormal covector from this formula or alternatively just like in the Euclidean case: $B_k=\epsilon_{kij}U^iN^j$ but we have to raise the index to a vector by changing the sign of the $0$ component: $B^i=\eta^{ij}B_j$. The result is
\beq
  \vec{\mathbf{B}} =  \gamma \langle -\sin(\gamma\omega \tau), \cos(\gamma\omega \tau), -v \rangle
                    = \frac{1}{\sqrt{c^2-a^2}}\langle -c\sin(\omega\omega \tau), c\cos(\gamma\omega \tau), -a \rangle \,,
\eeq
with torsion, now called the Fermi rotation of the spatial frame, equal to
\beq
 \Omega = \gamma^2\omega  = \frac{c}{c^2-a^2}\,.
\eeq

Since we assume $c>0$, the Fermi rotation or ``spin angular velocity" is positive, which means the two normal vectors rotate in the positive counter-clockwise azimuthal direction with respect to Fermi-Walker transported gyro axes, so the latter rotate backwards (retrograde direction) with respect to the Frenet-Serret spatial frame. We use the terminology spin angular velocity since this describes the precession of the direction of the electron spin in a circular orbit around the nucleus.

The fact that the proper orbital and spin angular velocities are not the same but instead related by the gamma factor
\beq
        \Omega = \gamma (\gamma \omega)
\eeq
means that they are not synchronized, so that while the Frenet-Serret spatial vectors return to the same orientation after each revolution of the helix, the Fermi-Walker transported vectors do not, but fall behind (retrograde direction) the Frenet-Serret  vectors by the difference frequency. 
\beq
    \Omega_{\rm Thomas} = \Omega-\gamma\omega = (\gamma-1) \gamma\omega\,.
\eeq
This is the famous Thomas precession of the electron spin in a circular orbit around the nucleus. The total angle per revolution is this frequency times the proper orbital period $2\pi/(\gamma\omega)$, namely
\beq
     \Delta\phi_{\rm thomas} = 2\pi (\gamma-1)
         = 2\pi \left( \frac{c}{c^2-a^2}\right)\overset{a/c \ll 1}{\rightarrow} \pi \frac{a^2}{v^2} \,.
\eeq
If we extend the two normal vectors to a coordinate system on the local rest space they span, we get the following representation of nearby points
\beq
   \Vec{x} =\Vec{x}(\tau) + X_{\rm FS} \Vec{N} + Y_{\rm FS} \Vec{N} \,,
\eeq
in terms of which one can parametrize a circle of radius $b$ about the spatial origin in the local rest space by the usual polar angle
\beq
    (X_{\rm FS} ,Y_{\rm FS} ) = (b \cos(v),b\sin(v)) \,.
\eeq
Our Lorentz corkscrew surface is just the world sheet of this circle.

Before leaving this analysis of the helical world line, it is worth finishing it off by introducing its Fermi coordinate system. We just have to reverse the Fermi angular velocity of the frame vectors in the local rest space to get a Fermi-Walker propagated orthonormal frame along the world line. We can then introduce spatial Fermi coordinates $(X,Y)$ based on the new spatial frame and use the proper time $\tau$ for the new time coordinate. First we rotate in the opposite direction by the angular velocity $\Omega$ to obtain a Fermi-Walker transported spatial frame 
\beq
  \begin{pmatrix} \Vec{N}_{\rm F} &  \Vec{B}_{\rm F} \end{pmatrix}
 =   \begin{pmatrix} \Vec{N} &  \Vec{B} \end{pmatrix}
       \begin{pmatrix} \cos(-\Omega\tau) &   \sin(-\Omega\tau) \\  -\sin(-\Omega\tau) &  \cos(-\Omega\tau) \end{pmatrix}
\eeq
and then we use this to locate points off the helical world line at $\Vec{x}(\tau)$ in the corresponding local rest space  with orthonormal coordinates in that local rest space with respect to the new spatial frame
\beq
   \Vec{x} =\Vec{x}(\tau) + X \Vec{N}_{\rm F} + Y \Vec{N}_{\rm F} \,.
\eeq
This is the coordinate transformation giving the inertial coordinates $(x,y,t)$ in terms of the Fermi coordinates $(X,Y,\tau)$.
One can easily visualize this through an animation along the curve of the local rest space with the spatial Fermi grid using a computer algebra system and 3d graphics. This shows clearly the Thomas precession very nicely \cite{bobFermipage}.
Of course realistic relativistic effects are usually very small, so only by exaggerating them enormously can one ``see them."

\section{World tube of a circle in the local rest space of the helix}

We can  boost a horizontal circular plane curve to the local rest space along this timelike helical central world line leading to a helical tube given by the parametrized surface by specializing the hyperbolic tilt parameter $\iota$ to be the rapidity $\beta$
\begin{eqnarray}
  \begin{pmatrix} x\\ y\\ t \end{pmatrix}
  &=&  (a+b\cos(v))\, \hat{\boldsymbol{\rho}}
     +    b\sin(v) \left( \cosh(\iota)\, \hat{\boldsymbol{\phi}} +\sinh(\iota) \,\hat{\mathbf{t}}\right)
     + c u\, \hat{\mathbf{t}}
\nonumber\\
  &=&  (a+b\cos(v)) 
  \begin{pmatrix} \cos(u)\\ \sin(u)\\ 0  \end{pmatrix}
\nonumber\\
  &&\qquad  + b\sin(v)   \left( \cosh(\iota) \begin{pmatrix} -\sin(u)\\ \cos(u)\\ 0  \end{pmatrix}
                      +\sinh(\iota) \begin{pmatrix} 0\\ 0\\ 1  \end{pmatrix} \right)
 + c u \begin{pmatrix} 0\\ 0\\ 1  \end{pmatrix}
\nonumber\\
 &=&   \begin{pmatrix} \left(a+b\cos(v)\right) \cos(u) -b\cosh(\iota)\sin(v)\sin(u)\\
                        \left(a+b\cos(v)\right) \sin(u) +b\cosh(\iota)\sin(v)\cos(u)\\
                        \phantom{\left(a+b\cos(v)\right) \sin(u)+} \kern-3pt b\sinh(\iota)\sin(v)
                         + c u  \end{pmatrix}
\nonumber\\
 &=&   \begin{pmatrix}\displaystyle \left(a+b\cos(v)\right) \cos(u) -\frac{bc}{\sqrt{c^2-a^2}}\sin(v)\sin(u)\\
                      \displaystyle   \left(a+b\cos(v)\right) \sin(u) +\frac{bc}{\sqrt{c^2-a^2}}\sin(v)\cos(u)\\
                       \displaystyle  \phantom{\left(a+b\cos(v)\right) \sin(u)+} \kern-3pt\frac{ab}{\sqrt{c^2-a^2}}\sin(v)
                         + c u  \end{pmatrix}
\,.
\end{eqnarray}
This is the relativistic cavatappo 2.0 surface, new and improved! To have a concrete example for illustrations we can a unit circular cross-section radius $b=1$ and set the speed parameter to be 1/2 the speed of light: $v=1/2$, $\gamma=2/\sqrt{3}, \beta={\rm arctanh}(1/2)$ which is accomplished by setting  $a=2b=2$, $c=4b=4$. The Euclidean angle of inclination is then $\eta=\arctan(2)\approx 63.4$ degrees. Because of the timelike condition, this angle is necessarily greater than 45 degrees so the helical tube is stretched out in the vertical direction compared to the Euclidean cavatappo surface. We could also consider a spacelike helix for the central curve, which would be more similar in shape to the canonical orthogonal cavatappo 2.0 surface illustrated above.  

Note that the $v$ coordinate lines on the surface are the result of intersecting the surface with the local rest spaces of this construction, while the $u$ coordinate lines are tied to the Fermi rotation of the two normal vectors, thus giving a nice visualization of the Frenet-Serret frame. If one introduces a new parametrization of the surface using the Fermi coordinate system, keeping the $u$ coordinate but translating $v$ to follow the Fermi-Walker propagated axes, one sees in the residual rotation of the new $u$ coordinate lines the remaining Thomas precession of the local rest space. The new grid is a concrete visual manifestation of the Fermi coordinates and the Thomas precession.

\begin{figure}
\typeout{*** EPS figure cavatappo relativistic quarter turn}
\begin{center}
%\vglue-2cm
 \includegraphics[scale=0.30]{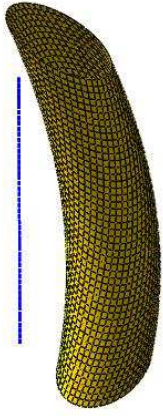}
\end{center}
%\vglue-2cm

\caption{One quarter turn of the Lorentz cavatappo together with the symmetry axis.
}\label{fig:cavatappolorentz}
\end{figure}

Evaluating the metric on the surface leads to a simple modification of the Euclidean result $c^2\to -c^2,v\to-v$ if one reinterprets $\sqrt{c^2+a^2}$ as $\sqrt{|c^2+a^2|}$
\begin{eqnarray}
ds^2&=&-\frac{ (c^2-a^2 -ab\cos(v))^2 -b^2 c^2}{c^2-a^2}du^2\nonumber\\
  &&\qquad    + \frac{2b^2c}{\sqrt{c^2-a^2}}  \, du \, dv + b^2 dv^2
\,.
\end{eqnarray}

It is helpful to introduce some decomposition quantities for both cases, with appropriate notation for each.
Adapting the decomposition to be orthogonal with respect to the symmetry group orbits, namely the $u$ coordinate lines or  ``threads," gives the threading decomposition (complete the square on $du$)
\begin{eqnarray} \label{eq:metricML}
 ds^2
&=& g_{uu} \, du^2 +2 g_{uv} \, du\, dv + g_{vv}\, dv^2\\
&=& -M^2 (du-M_v dv)^2 +\gamma_{vv} \, dv^2
\,,\nonumber\\
&=& -(\omega^{\top})^2 +(\omega^{\hat v})^2 \,,
\end{eqnarray}
which requires computation of the following quantities:
\begin{eqnarray}
  M &=&(-g_{uu})^{1/2}
=  \left(\frac{\strut (c^2-a^2-ab\cos(v))^2 -b^2c^2}{c^2-a^2}\right)^{1/2}
\,,
\nonumber\\     
M_v &=& \frac{g_{uv}}{M^2} =  \frac{b^2c}{M^2\sqrt{c^2-a^2}} 
\,,
\qquad
   M^v = M_v/g_{vv}= \frac{c}{M^2\sqrt{c^2-a^2}}
\,,
\nonumber\\
  \gamma_{vv} &=& g_{vv} -g_{uv}^2/M^2 
  =
%\frac{b^2 (a^2+c^2+ab\cos(v))^2}{\left((a^2+c^2)(a+b\cos(v))^2+c^2+b^2c^2\sin^2(v) \right)}
\frac{b^2 (c^2-a^2-ab\cos(v))^2}{\strut (c^2-a^2-ab\cos(v))^2 -b^2c^2}
\,,
\end{eqnarray}
where the metric determinant is
\beq
   -\det(g_{ij}) =  M^2 \gamma_{vv} = \frac{b^2 (c^2-a^2-ab\cos(v))^2}{\strut (c^2-a^2)} 
\,.
\eeq
This corresponds to the orthonormal frame and dual frame
\begin{eqnarray}
&& \Vec{e}_{\top} = M^{-1/2} \partial_u\,,\
   \Vec{\epsilon}_{\hat v} = \gamma_{vv}^{-1/2} (\partial_v+M^v\partial_u)\,,
\nonumber\\   
&& \omega^{\top} = M (du-M_v\, dv)\,,\
   \omega^{\hat v} = \gamma_{vv}^{1/2} \omega^v\,.
\end{eqnarray}
Fig.~\ref{fig:threadingsplit} illustrates this orthogonal decomposition of the tangent space.

The differential surface volume element is
$$
  dS = |\det(g)|^{1/2} du\, dv 
     = \frac{b (c^2-a^2-ab\cos(v))}{\sqrt{c^2-a^2}}\,,
$$
and the total surface area of one revolution of the surface is
$$
           S = (2\pi b) (2\pi\sqrt{c^2-a^2}) = (2\pi b) T_o
\,.
$$
This is just the circumference $C=2\pi b$ of the cross-sectional circle times the total proper time along the central helix during one revolution (the proper period $T_o$), extending the theorem of Pappus to the Lorentzian case. The fact that the integral of the differential of surface area is trivial to integrate by hand shows how much simpler the geometry of the orthogonally tilted cavatappo surface is compared to the general case, which allows exact integration but which results in a formula full of elliptic functions that covers many computer screens.

One easily calculates the Riemann scalar curvature $R$ and the Gaussian curvature $K$ with a computer algebra system to be
\beq
 R =-2K= \frac{2a\cos(v) }{b(c^2-a^2-ab\cos(v))} \,.
\eeq
which has the same sign as $\cos(v)$, namely positive on the outer hemisurface and negative on the inner hemisurface, while vanishing at the polar helices. The Gaussian curvature has the opposite sign due to the Lorentzian signature change in the classical Gauss equation for a hypersurface in a Riemannian manifold.  (See Eq.~(21.75) of Misner, Thorne and Wheeler \cite{mtw},  doubly contracted to the curvature scalar and applied to a flat enveloping spacetime.)

% cavatappo20-threadingsplit

\begin{figure}
\typeout{*** EPS figure cavatappo threading splitting}
\begin{center}
%\vglue-2cm
\includegraphics[scale=0.35]{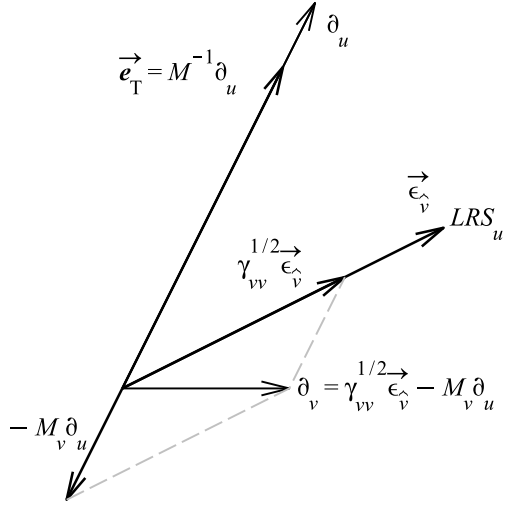}
\end{center}
%\vglue-2cm

\caption{The orthogonal decomposition of the tangent space with respect to the time lines, showing the spacelike meridian direction $\partial_v$ as horizontal for convenience. 
}\label{fig:threadingsplitLorentz}
\end{figure}

The slicing decomposition which is orthogonal with respect to the circular meridians (complete the square on $dv$) is
\begin{eqnarray}\label{eq:metricNL}
 ds^2 &=& -N^2 du^2 +g_{vv} \, (dv+N^v du)^2
\,,\nonumber\\
&=&- (\theta^\bot)^2 +(\theta^{\hat v})^2 \,,
\end{eqnarray}
which requires computation of the following quantities:
\begin{eqnarray}
  N &=&(|g_{uu}-g_{uv}^2/g_{vv}|)^{1/2} 
%= \left( (a+b\cos(v))^2-\frac{b^2 c^2 \cos^2(v)}{a^2+c^2} +c^2\right)^{1/2} 
=  \frac{c^2-a^2 -a b \cos(v)}{\sqrt{c^2-a^2}} 
\nonumber\\
    N^v &=& \frac{g_{uv}}{g_{vv}} =   \frac{c}{\sqrt{c^2-a^2}}
\,,\ 
N_v = g_{uv}=g_{vv} N^v=\frac{b^2 c}{\sqrt{c^2-a^2}}
\,,
\nonumber\\
g_{vv} &=& b^2
\,,
\end{eqnarray}
while the metric determinant is also expressable as
\beq
   -\det(g_{ij})
    = N^2 g_{vv} = M^2 \gamma_{vv}  \frac{b^2 (c^2-a^2-ab\cos(v))^2}{\strut (c^2-a^2)}
\,.
\eeq

\begin{figure}
\typeout{*** EPS figure cavatappo threading splitting}
\begin{center}
%\vglue-2cm
\includegraphics[scale=0.35]{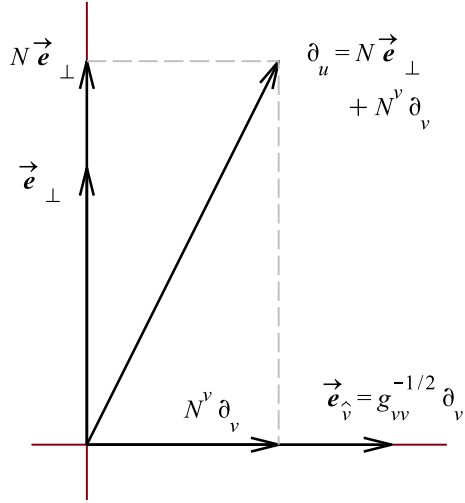}
\end{center}
%\vglue-2cm

\caption{The orthogonal decomposition of the tangent space with respect to the spacelike meridians, showing the meridian direction $\partial_v$ as horizontal for convenience. 
}\label{fig:slicingsplit}
\end{figure}

This orthogonal decomposition of the metric corresponds to the orthonormal frame and dual frame
\begin{eqnarray}
&& \Vec{e}_{\bot} = N^{-1/2} (\partial_u-N^v\partial_v)\,,\
  \Vec{e}_{\hat v} = g_{vv}^{-1/2} \partial_v\,,
\nonumber\\   
&& \omega^{\top} = N \, du \,,\
   \omega^{\hat v} = g_{vv}^{1/2} (dv+N^v du)\,.
\end{eqnarray}
Fig.~\ref{fig:slicingsplit} illustrates this orthogonal decomposition of the tangent space.

% cavatappo20-slicingsplit

\begin{figure}
\typeout{*** EPS figure cavatappo Fermi grid}
\begin{center}
%\vglue-2cm
\includegraphics[scale=0.35]{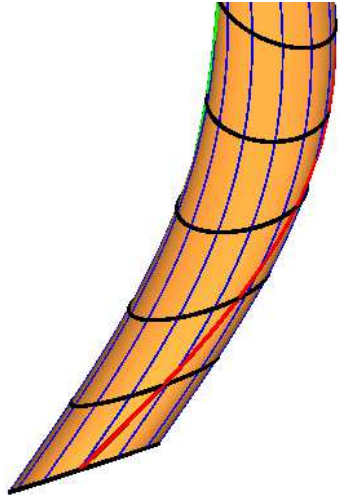}
\end{center}
%\vglue-2cm

\caption{The coordinate grid induced from the Fermi coordinate grid along the central helix in the 3-dimensional Minkowski spacetime, with the same spatial meridians but new timelike parallels orthogonal to those meridians.
}\label{fig:fermigrid}
\end{figure}

The new time lines orthogonal to the spatial meridians are simple
\begin{eqnarray}
  g_{\hat v\hat v}^{-1/2}\theta^{\hat v} &=& dv+N^v du =d\left(v-\frac{c}{\sqrt{c^2-a^2}} u\right) \equiv
  dv_\bot
\nonumber\\
  &\to& v_\bot =v-\frac{cu}{\sqrt{c^2-a^2}}\,.
\end{eqnarray}
The new grid associated with the orthogonal coordinates $(u,v_\bot)$ has the same spatial meridians but new time lines associated with Fermi coordinates on the 3-dimensional Minkowski spacetime defined near the central helix on which they are based \cite{fermicircles}. The residual rotation of this ``nonrotating" grid compared to the spacetime inertial time lines after each successive revolution of the helix is a visual representation of the Thomas precession which rotates the spin of the electron in the plane of its circular orbit. Fig.~\ref{fig:fermigrid} illustrates this new coordinate grid.

This orthogonal grid is in fact formed from the two families of lines of curvature of the surface. Evaluating the extrinsic curvature, one finds the principal curvatures are
$$
  k_1 = \frac1b\,,\ k_2 = -\frac{a\cos(v)}{c^2-a^2-ab\cos(\varphi)}
$$
along the meridians and their orthogonal trajectories respectively, and
related to the scalar curvature by $ R=-2K = -2(k_1k_2)^{-1}$. Note the additional minus sign that occurs comparing twice the Gaussian curvature to the intrinsic Riemann scalar curvature when dealing with a spacelike hypersurface in a Lorentzian spacetime.

The zero screw-angular momentum curves are the spatial curves orthogonal to the existing parallels or time lines. These are instead defined by
\begin{eqnarray}
   \omega^{\top} &=& M (du-M_v\, dv)=0 \to \frac{du}{dv}= M_v 
 = \frac{b^2c}{M^2\sqrt{c^2-a^2}} 
\nonumber\\
  &\to& u=\int_{v_0}^v 
   \frac{\sqrt{c^2-a^2}}{(c^2-a^2-ab\cos(\varphi))^2 -b^2c^2}
 d\varphi
 \,,
\end{eqnarray}
which can be integrated exactly but one has to adjust the periodic antiderivative to be cumulative by adding the total increment $\Delta u>0$ for one revolution in $v$ after each such revolution. These curves are unbounded, wrapping around the tube forever. These are spatial geodesics. The increment $\Delta u$ is the synchronization gap \cite{clock}, nonzero since the time lines do not admit closed orthogonal curves which would correspond to a global synchronization of the coordinate time $u$.

If $A(v)$ is the periodic antiderivative on $[-\pi,\pi]$ for which $A(0)=0$, then the corresponding cumulative antiderivative is
\begin{eqnarray}
   {\rm ShiftTo0Pi}(v)&=& v-2\pi \,{\rm floor}\left(\frac{v}{2\pi}\right)\,,
\nonumber\\
   AC(v) &=& A({\rm ShiftTo0Pi}(v))+ 2 \Delta u \,{\rm floor}(\left( \frac{v+\pi}{2\pi}\right) \,.
\end{eqnarray}

\section{Lorentz Geodesics}

The geodesic equations are also very closely related to the previous Euclidean case by the transformation $c^2\to-c^2, v\to-v$
\begin{eqnarray}
&&\mathcal{D} \frac{d^2 u}{d \tau^2}
 +\frac{a b c\sin(v) }
       {\sqrt{c ^{2}-a ^{2}}} \left(\frac{d u}{d t } \right)^2
    +2 a b \sin \left(v \right)
     \frac{d u}{d t }  \frac{d v}{d t }
=0
\,,\nonumber\\
&&\mathcal{D} \frac{d^2 v}{d \tau^2}
 -\frac{a \sin \left(v \right) \left(\left(c^2-a^2-a b \cos \left(v \right)\right)^{2}-b^2c^2\right)}
          {b \left(c^2-a^2\right)} 
     \left(\frac{d u}{d t } \right)^2
\nonumber\\
&&\quad +\frac{2 a b c \sin \left(v \right)}{\sqrt{c^2-a^2}} 
      \frac{d u}{d t }  \frac{d v}{d t }
=0 \,,
\end{eqnarray}
where
\beq
\mathcal{D} = c^2-a^2-ab\cos(v) \,.
\eeq
The meridians $u=u_0$ are again obvious closed circular spatial geodesics.

\section{Constants of the motion, initial data problem}

If $(u(\lambda),v(\lambda))$ is an affinely parametrized geodesic of this metric on the surface, the same formulas hold for the conserved screw-angular momentum threading quantities already discussed for the Euclidean case but with sign changes, and now the geodesics are also characterized by being either timelike, spacelike  or lightlike (null).

\begin{figure}[p]
\typeout{*** EPS figure cavatappo}
\begin{center}
%\vglue-2cm
\includegraphics[scale=0.45]{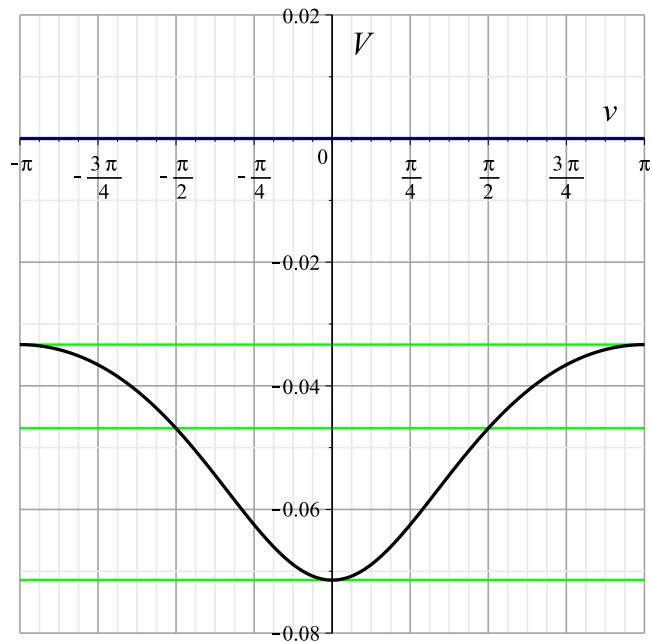}
\end{center}
%\vglue-2cm

\caption{The Lorentz screw-centrifugal potential. Shown in green are the energy levels of the inner and outer equators and of the northern and southern polar helices, all timelike. Positive energy levels $E>0$ correspond to spacelike geodesics, $E=0$ to null, and negative $E<0$ to timelike geodesics. Infinite energy corresponds to the zero screw-angular momentum limit. All motion with energy above the inner equator value is unbounded in $v$. 
}\label{fig:screw-centrifugalpot-lorentz}
\end{figure}

The tangent vector in the threading decomposition satisfies
\begin{eqnarray}
  U^{\hat u} = M \left(  \frac{du}{d\lambda} - M_v \frac{dv}{d\lambda} \right ) 
\quad 
  U^{\hat v} = \gamma_{vv}^{1/2}  \frac{dv}{d\lambda} \,,
\end{eqnarray}
and
\begin{eqnarray}  \label{eq:velocityrelation-Lorentz}
 U^u = \frac{du}{d\lambda} = M^{-1} U^{\hat u} + M_v \gamma_{vv}^{-1/2}   U^{\hat v}\,,
\quad %\nonumber\\
  U^v= \frac{dv}{d\lambda} = \gamma_{vv}^{-1/2}  U^{\hat v}
\,. 
\end{eqnarray}
The conserved screw-angular momentum is 
\beq
  \ell  =U_u = g_{uu}U^u+g_{uv}U^v
   = M U^{\hat u} = M^2 \left( \frac{du}{d\lambda}-M_v \frac{dv}{d\lambda} \right)
\,, 
\eeq
while the conserved energy $E >0,=0,<0$ (spacelike, null, timelike) is 
\beq
  - {\frac12}\left( U^{\hat u} \right)^2 +   {\frac12}\left( U^{\hat v} \right)^2
 =  {\frac12}\left( U^{\hat v} \right)^2 -  {\frac12}\frac{\ell ^2}{M^2} =  E 
\eeq
with the effective potential for the motion along the meridians
\beq
   V =  -{\frac12}\frac{\ell ^2}{M^2}
=-\frac{(c^2-a^2)\ell^2}{2(c^2-a^2-ab\cos(v))^2 -b^2c^2} \,.
\eeq
Fig.~\ref{fig:screw-centrifugalpot-lorentz}
shows the nonzero screw-angular momentum potential for the orthogonally tilted Lorentz cavatappo 2.0 pasta surface for the canonical parameter values.

\begin{figure}
\typeout{*** EPS figure cavatappo timelike initial data}
\begin{center}
%\vglue-2cm
\includegraphics[scale=0.35]{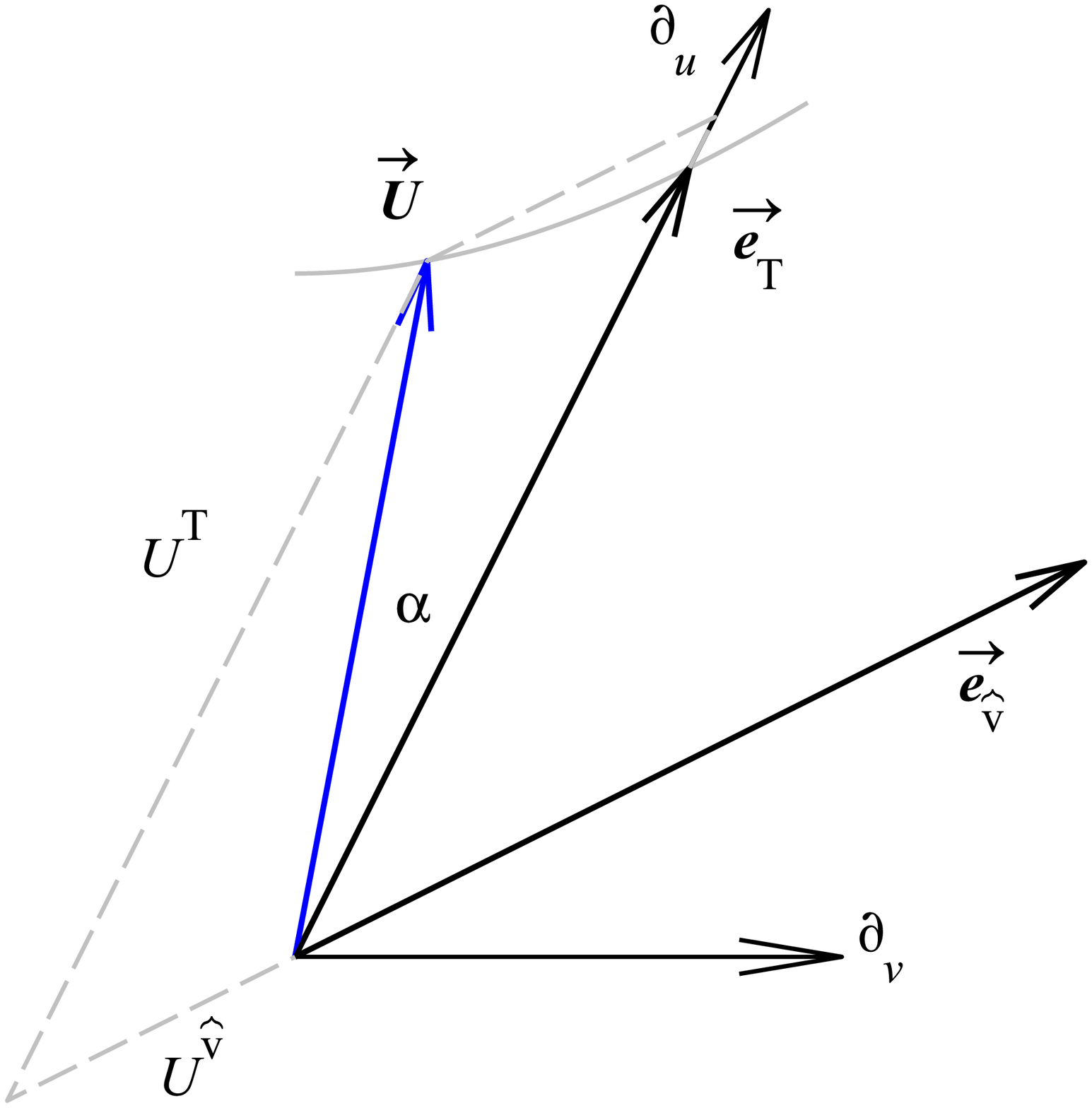}
\end{center}
%\vglue-2cm

\caption{Orthogonal decomposition of initial data in the vertical tangent plane at the outer equator for a unit timelike tangent vector $\Vec{U}$ with respect to the orthonormal frame adapted to the time lines along $\partial_u$, showing the $v$ coordinate direction $\partial_v$ horizontal only for convenience.
}\label{fig:initialdatatimelike}
\end{figure}

The initial data problem here is complicated by the three classes of geodesics: $\vec{\mathbf{U}}\cdot \vec{\mathbf{U}} = \epsilon |2E|, \epsilon={\rm sgn}(E)=1,0,-1$ (respectively spacelike, null and timelike). Now we must use either a hyperbolic angle with respect to the timelike direction $\Vec{e}_{\top}$ for timelike geodesics, or with respect to the spacelike direction $\Vec{e}_{\hat u}$, while null geodesics are a special case $|U^{\top}|=|U^{\hat u}|$. For a nonnull geodesic the unit tangent is obtained by dividing the tangent $\vec{\mathbf{U}}$ by its length $2E$. Fig.~\ref{fig:initialdatatimelike} shows the orthogonal decomposition of the initial data for the case of a timelike geodesic.

%cavatappo20-initialdata

The spacelike case is similar to the Euclidean case so consider the timelike case, where this normalized tangent can be interpreted as a ``4-velocity" in the 2-dimensional spacetime.
Initial data can be posed at the origin of coordinates on the outer equator at the prime meridian
\beq
 \vec{\mathbf{U}}   = U^{\hat u} \Vec{e}_{\top} + U^{\hat v} \Vec{\epsilon}_{\hat v}\,,\quad
 (U^{\hat u} , U^{\hat v}) =  (2E)(\cosh\alpha,\sinh\alpha) \,.
\eeq
The coordinate velocities are then
\begin{eqnarray}  \label{eq:velocityrelationinitialLorentz}
\frac{du}{d\lambda} = 2E\left( M^{-1} \cosh\alpha - M_v \gamma_{vv}^{-1/2}   \sinh\alpha\right)\,,
\quad
\frac{dv}{d\lambda} = \gamma_{vv}^{-1/2}  \sinh\alpha
\,, 
\end{eqnarray}
in terms of which the conserved screw-angular momentum is
\beq
  \ell = M (2E) \cosh\alpha \,.
\eeq
For numerical solution a unit arclength parametrization is simplest.

\begin{figure}[p]
\typeout{*** EPS figure cavatappo null geos}
\begin{center}
%\vglue-2cm
\includegraphics[scale=0.45]{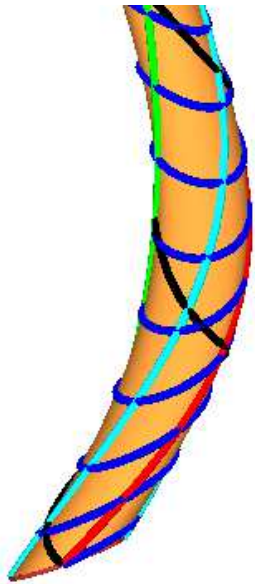}
\end{center}
%\vglue-2cm

\caption{The two null geodesics leaving the origin of coordinates at intersection of the prime meridian (orange) and the outer equator (red). The prograde orbit (blue) is much slower to return to the outer equator compared to the retrograde orbit (black).
}\label{fig:nullgeos}
\end{figure}

The single null geodesic at the origin of coordinates corresponds to the initial data $(U^{\hat u} , U^{\hat v})=(1,\pm 1)$ so the initial data can be chosen so that 
\begin{eqnarray}  \label{eq:velocityrelationinitialnull}
\frac{du}{d\lambda} &=& \left( M^{-1}  \mp M_v \gamma_{vv}^{-1/2} \right)\,,
\quad
\frac{dv}{d\lambda} = \pm\gamma_{vv}^{-1/2}
\,, 
\qquad
  \ell = M \,.
\end{eqnarray}
Fig.~\ref{fig:nullgeos} shows a pair of null geodesics leaving the origin of coordinates in opposite directions. The prograde orbit leaves with $U^{\hat v}=-1$ and is much slower in returning to the outer equator compared to the retrograde orbit with $U^{\hat v}=-1$, which is tightly wound around the surface. This is an exaggerated illustration of an effect similar to the Sagnac effect in which two photons in circular orbit around a rotating black hole return to the same spatial point fixed in the stationary coordinate grid at different times \cite{clock}.

The tangent plane to the surface at the origin is vertical. One can also determine the initial data there which is aligned with the inertial coordinates $(t,y)$ (vertical, horizontal) and examine the corresponding timelike and spacelike geodesics. For example, the condition
\beq
 \frac{U^v}{U^u} = - \frac{(a+b)\sqrt{c^2-a^2}}{bc} \,.
\eeq
determines the boost angle of $\partial_u\sim {\partial R}/{\partial u}(0,0)$ back to the vertical direction $\langle 0,0,1\rangle$, while
\beq
 \frac{U^v}{U^u} = - \frac{c\sqrt{c^2-a^2}}{ab} \,.
\eeq
determines the boost angle of $\partial_v\sim ({\partial \Vec{x}}/{\partial v})(0,0)$ back to the horizontal, where $\sim$ indicates that these tangent vectors agree in the enveloping geometry of the whole space. Using the respective boost angles determined from these relations in the initial data relations determines the corresponding initial data for the numerical geodesic machine. 

\section{Concluding remarks}

Toy examples are very useful for making mathematical and physical concepts come alive and provide mental pictures that aid our understanding. While the classical theory of surfaces may be interesting to some as it has been presented for nearly a century, for those of us also interested in theoretical physics at a technical or even pedestrian level, finding unorthodox new examples like this one, whether it be pasta or a relativistic circularly orbiting ball, is a real joy. Computer algebra systems make the necessary calculations reasonable, and generate impressive graphics to drive home the consequences of those calculations. I hope a few others take some inspiration from this discussion.

%%%%%%%%%%%%%%%%%%%%%%%%%%%%%%%%%%%%%%%%%%%%%%%%%%%%%%%%%%%%%%%%%%%%%%%%%%%%%%%%%%%%%%%%%%%%%%%%%%%%%%%%%%%%%%%%%
%\newpage
\section*{Acknowledgements}

This work would not have been possible without the initial stimulus from my friend and colleague Klaus Volpert and my long history of ignoring the classical differential geometry of curves and surfaces while playing games in mathematical general relativity. Thanks also to Chris Rorres for his inspiring 2012 talk on the Archimedes screw which sparked the light bulb moment leading to the cavatappo 2.0 generalization of earlier work. Finally without a computer algebra system, this work would never have even begun. All the figures were produced and the underlying calculations were performed with Maple (TM) \cite{maple}.

\end{document}